\documentclass[reqno,11pt,a4paper]{amsart}

\usepackage{amsmath}
\usepackage{amsthm}
\usepackage{amsfonts,amssymb}
\usepackage{bbm}
\usepackage{multirow}
\usepackage{latexsym}
\usepackage{mathrsfs}
\usepackage[all]{xy}
\usepackage{url}

\usepackage{float}
\usepackage{hyperref}
\usepackage{amssymb,longtable}
\usepackage{fontenc}
\usepackage{array,etoolbox}
\preto\tabular{\setcounter{magicrownumbers}{0}}
\newcounter{magicrownumbers}
\newcommand\rownumber{\stepcounter{magicrownumbers}\arabic{magicrownumbers}}
\usepackage{pdfsync}

\usepackage{amsthm}
\usepackage{hyperref}
\usepackage{enumerate}
\usepackage{url}

\usepackage{longtable}
\usepackage{booktabs}
\usepackage{colortbl}
\newcommand{\evnrow}{\rowcolor[gray]{0.95}}
\newcommand{\oddrow}{}

\usepackage[margin=2.51cm]{geometry}

\usepackage{pdflscape,fullpage}
\usepackage{longtable}
\usepackage{booktabs}
\usepackage{color}
\usepackage{arydshln}
\usepackage{verbatim}
\usepackage{enumitem}

\usepackage{mathdots}

\usepackage[english]{babel}
\usepackage{lipsum}
\usepackage{amsthm}
\usepackage{todonotes}
\usepackage[vcentermath]{youngtab}

\theoremstyle{plain}
\newtheorem{lema}{Lemma}[section]

\newtheorem{thrm}[lema]{Theorem}

\theoremstyle{remark}

\theoremstyle{definition}
\newtheorem{dfn}[lema]{Definition}
\newtheorem{eg}[lema]{Example}

\newcommand{\om}{\omega}

\def\ZZ{{\mathbb Z}}
\def\QQ{{\mathbb Q}}
\def\CC{{\mathbb C}}

\def\AA{{\mathbb A}}

\def\PP{{\mathbb P}}
\def\wgr{\mathrm{wGr(2,5)}}
\def\Gr{\mathrm{Gr}}
\def\wp2{\mathrm{w}($\PP^2\times \PP^2$)}

\def\Spec{{\mathrm{Spec}}}
\def\Proj{\mathrm{Proj}}

\def\sm{\mathrm{sm}}

\def\w{{\mathrm{w}}}

\def\Pf{{\mathrm{Pf}}}

\def\orb{\mathrm{orb}}
\def\GCD{\mathrm{GCD}}
\def\Max{\mathrm{maximum}}

 \newcommand{\into}{\hookrightarrow}

 \newcommand{\Oh}{\mathcal O}

 \newcommand{\sB}{\mathcal B}

 \newcommand{\Si}{\Sigma}

\newcommand{\PxP}{\PP^2 \times \PP^2}

\begin{document}

\author[M.I.~Qureshi]{Muhammad Imran Qureshi}
\address{Muhammad Imran Qureshi,
Department of Mathematics, King Fahd University of Petroleum \& Minerals (KFUPM), Dhahran
31261, KSA }
\email{imran.qureshi@kfupm.edu.sa}
\title{Polarized rigid del Pezzo surfaces  in low codimension }


\begin{abstract}
We provide explicit graded constructions of orbifold del Pezzo surfaces with rigid orbifold points of type $\left\{k_i\times\frac{1}{r_i}(1,a_i): 3\le r_i \le 10,k_i \in \ZZ_{\ge 0}\right\}$; as well-formed and quasismooth varieties embedded in  some weighted projective space. In particular, we present a  collection of 147 such surfaces such that their image under 
their anti-canonical embeddings can be described by using one of the following sets of equations:
a single equation, two linearly independent equations, five maximal Pfaffians of \(5\times 5\) skew symmetric matrix, and nine \(2\times 2\) minors of size 3 square matrix. This is a complete classification of such surfaces under certain carefully chosen bounds on the weights of ambient weighted projective spaces and it is largely based on detailed computer-assisted searches by using the computer algebra system \textsc{magma}.    \end{abstract} 
\maketitle

\section{Introduction}

A del Pezzo surface is  a two dimensional algebraic variety with an ample anti-canonical divisor class. The classification of nonsingular  del Pezzo surfaces is  well known and  there are ten deformation families of such   surfaces: \(\PP^1 \times \PP^1\), \(\PP^2\) and the  blow up of \(\PP^2\) in \(d\) general points for \(1 \le d\le 8\).
An  orbifold   del Pezzo surface $X$ is a del Pezzo surface with at worst isolated    orbifold points, classically  known as a log del Pezzo surface with cyclic quotient singularities. We describe  \(X\) to be locally qGorenstein($qG$)-rigid if it contains only rigid isolated orbifold points, i.e. the orbifold points are rigid  under \(qG\)-deformations.  If it admits a \(qG\)-degeneration to a normal toric del Pezzo surface then it is called a del Pezzo surface of class TG. The Fano index  of \(X\) is the largest integer \(I \) such that \(K_X=ID\) for an element \(D\) in the class group of \(X\).

The classification of orbifold del Pezzo surfaces has been an interesting area of research from various points of view, such as the existence of Kahler--Einstein metric \cite{dp-exceptional,dp-zoo}.    
Recently the classification of orbifold del Pezzo surfaces has received much attention, primarily due to the mirror symmetry program for Fano varieties  by Coates, Corti et al
\cite{CCGGK}. The mirror symmetry   for orbifold del Pezzo surface  has been  formulated  in   \cite{dp-AMS} in the form of a conjecture expecting a one to one correspondence between mutation equivalence  classes of Fano polygons with the ($qG$)-deformation equivalence  classes of  locally \(qG\)-rigid del Pezzo surfaces of class TG.   Therefore the construction  of rigid orbifold del Pezzo surfaces has important links with the mirror symmetry due to this conjecture. The conjecture has been proved  for smooth del Pezzo surfaces by Kasprzyk, Nill and Prince in~\cite{dp-KNP}.  
    Corti and Heuberger ~\cite{dp-CH} gave the classification of locally \(qG\)-rigid del Pezzo surfaces  with \(\frac13(1,1)\) singular points. The del Pezzo surfaces with a single orbifold point of type  \(\frac 1r(1,1)\) have been classified  by Cavey and Prince \cite{dp-CP}. The mutation equivalence classes of Fano polygons with rigid singularities of type        \begin{equation}\label{eq:basket_CK}\left\{k_1\times \frac13(1,1),k_2\times \frac16(1,1): k_1> 0,k_2 \ge 0 \right\} \text{ and } \left\{k\times \frac15(1,1): k>0\right\}\end{equation}
have been computed in \cite{Fano-CK}. This  is equivalent to  the classification of del Pezzo surfaces of class TG with  the above given baskets;  though it may be missing  surfaces which do not admit a toric degeneration and  having one of the above type of baskets of singularities. By using birational techniques,  the classification of orbifold del Pezzo surfaces with basket consisting of a combination of \( \frac13(1,1)\) and \(\frac14(1,1)\) orbifold points was given by Miura \cite{dp-Miura}.

 In \cite{dp-CH} the  classification gave a total of   29 deformation families of del Pezzo surfaces with \(\frac13(1,1)\) orbifold points which were divided into  6 different cascades; one of the cascades was first studied by Reid and Suzuki in \cite{dp-RS}. Moreover,  good model  constructions for  all  29 surfaces were presented  as complete intersections inside the so called rep-quotient varieties (mainly simplicial toric varieties): A geometric quotient \(V//G\)  of a representation \(V \) of a complex Lie group \(G\) .   Among those, six of them can be described  as a hypersurface in \(\PP^3(a_i)\) or as a complete intersection in $\PP^4(a_i)$ or as complete intersection in weighted Grassmannian \(\w\Gr(2,5)\)\cite{wg}.
This motivated us to classify   rigid del Pezzo surfaces with      certain basket of  singularities  which can be described by relatively small sets of equations. 
\subsection{Summary of results}
We classify polarized rigid del Pezzo surfaces, under the bounds chosen in \ref{Bounds}, which contain  baskets of orbifold points   \[\left\{k_i\times \frac{1}{r_i}(1,a_i): 3\le r_i\le 10,  k_i \ge 0\right\};\]
such that their images under their anti-canonical embedding can be described by one of the following ways.
\begin{enumerate}
\item as a hypersurface, i.e. by a single weighted homogenous equation; \(X_d \into\PP^3(a_i)\).
\item as a codimension 2 weighted complete intersection, i.e. by 2 weighted homogeneous equations; \(X_{d_1,d_2}\into \PP^4(a_i).\) 
\item  as a codimension  3 variety described by using  five maximal Pfaffians of a \(5\times 5\) skew symmetric matrix;  $$X_{d_1,\ldots,d_5}\into \PP^5(a_i).$$ In other words they are weighted complete intersections in weighted Grassmannian \(\w\Gr(2,5)\) or (weighted) projective cone(s) over it \cite{wg,QS,QS-AHEP}. 
\item as a codimension 4 variety described by using  nine  \(2\times 2\) minors of a size 3 square matrix $$X_{d_1,\ldots,d_9}\into\PP^6(a_i).$$
Equivalently they are  weighted complete intersections in some  weighted \(\PxP\) variety or (weighted) projective cone(s) over it \cite{BKQ}.  

\end{enumerate}
 We summarize the classification in  form of the  following theorem. 
\begin{thrm}\label{MainResult} Let \(X\) be an orbifold del Pezzo surface having at worst a basket  $$\mathcal B=\left\{k_i \times \frac1r(1,a):3 \le r\le 10, k_i \ge 0\right\}$$of rigid orbifold points and  their image  \(X\into \PP(a_i)\) under their anti-canonical embedding can be  described as a hypersurface or as a codimension 2 complete intersection or  as a weighted complete intersection in \(\wgr\) or as a weighted complete intersection of weighted \(\PxP\) variety. Then, subject to  \ref{Bounds},  \(X\) is one of the del Pezzo surfaces listed in Tables ~\ref{tab-hypersurface}--\ref{tab-P2xP2}. In total there are 147  families of such del Pezzo surfaces, divided as follows in each codimension.

\[
\renewcommand*{\arraystretch}{1.2} \begin{tabular}{|c|c|c|c|}\hline
 Hypersurface & Complete intersection   & \(4 \times 4 \) Pfaffians & $2 \times 2$ Minors   
 \\\hline\hline
\evnrow  81 & 25 & 21 &20\\\hline

\end{tabular}\]
   \end{thrm}
\noindent    We construct these examples by first computing all possible candidate varieties with required basket of orbifold points using an  algorithmic approach  developed  in \cite{BKZ,QJSC}, under the   bounds given in \ref{Bounds}. In case of  codimension 1 and
2, the equations of these varieties are generic weighted homogeneous polynomials of given degrees. In cases of codimension 3 and 4 they are induced from the equations of the corresponding ambient weighted projective variety.  We perform a detailed singularity analysis of equations of these candidate varieties to prove the existence or non-existence of given candidate surface. We calculate the \(qG\)-deformation invariants like the anti-canonical degree \(-K_X^{2 }\) and first plurigenus \(h^0(-K_X)\) in all cases. We calculate their Euler number and  Picard rank in hypersurface case. In complete intersection case, we were able to calculate their Euler number and identify the non-prime examples, i.e. those with the  Picard  rank greater than 1 by computing their orbifold Euler number.      

The computer search used to find these surfaces, based on the algorithm approach of \cite{BKZ,QJSC} ,
is an infinite search. The search is usually performed in the order of increasing  sum of the weights \(W=\sum a_i\) of the ambient weighted projective spaces.  In each codimension and for each Fano index \(I\), we provide complete classification of rigid del Pezzo surfaces $X\subset \PP(a_i)$ satisfying $$W-I \le N  \text{ where } N \ge 50 .$$  If the last candidate example for computer search appears for \(W-I=q\) then we search for all cases with \(N=\Max(50,2q)\),  to minimize the possibility of any further examples. This indeed does not rule out a possibility of further other examples for larger value of \(W\) and \(I. \) Though It is  evident that for larger values of \(W\) most weights of \(\PP(a_i) \) will be larger than 10; highest local index of allowed orbifold points in our classification, consequently the basket of orbifold points will very likely contain orbifold points of local index \(r\ge 11\). In cases of hypersurfaces and complete intersections, the  classifications of tuples \((\underline {d_{j}};\underline {a_i})\) which
give rise to a quasismooth del Pezzo surfaces can be found in \cite{Mayanskiy}
and\cite{dp-Erik} ; where $\underline {d_{j}}$ denote the degrees of the defining equations and \(\underline {a_{i}}\) are weights of the ambient weighted projective space. These classifications of tuples can perhaps be analyzed to give the bound free proof of completeness of our results in codimension 1 and 2. Though their classification neither  contains computation of  any of the invariants like \(h^0(-K_X),-K_X^2\) and \(e(X)\) nor they  compute basket of orbifold points lying on those surfaces.

\subsection{Links with existing literature}
     A part of our search results recovers some   existing examples in the literature, though a significant subset   of them had not been  described earlier  in terms of equations. For example, the classification of  Fano polygons (equivalently of rigid del Pezzo surfaces of class TG) with basket of orbifold points \eqref{eq:basket_CK} is given in \cite{Fano-CK}. We give descriptions in terms of equations for six of their examples; listed as  \(14, 16,23, 85, 109 \text{ and } 130\) in our tables. We also  recover  the classical smooth del Pezzo surfaces of degrees \(1, 2,3,4,5,6 \text{ and }8\);   listed as \(3,2,1, 82,107,128\text{ and        }12 \) respectively in Tables \ref{tab-hypersurface}--\ref{tab-P2xP2}.  Moreover, 7 of the 29 examples from \cite{dp-CH} also appear in our list with one of them seemingly having  a new description as a complete intersection in a \(\w(\PxP)\) variety, listed as \(129\) in Table \ref{tab-P2xP2}.  Some  examples of Fano index \(1\) and \(2\) in codimension 3 and  4 given  in Table \ref{tab-Pf} and \ref{tab-P2xP2} can be found in  \cite{QMOC}, primarily appearing implicitly as a part of some infinite series  of orbifold del Pezzo surfaces.       
\section{Background and notational conventions}
\subsection*{Notation and Conventions}
\begin{itemize}
\item We work over the field of complex numbers \(\CC\). 
\item All of our varieties are projectively Gorenstein.
\item  For two  orbifold points where  \(\frac1r(1,a)=\frac1r(1,b)\)  we choose a presentation  \(\frac1r(1, \min (a ,b))\).  
\item In all the tables,  integers appearing as subscripts of \(X\) denote the degree of the defining equations of the given variety, where \(d^m\) means that there are \( m\) equations of degree \(d\).
Similarly,       \(\PP(\cdots,a_i^m,\cdots)\) means that there are \(m\) weights of degree \(a_i\).
\item We use the same notation for canonical divisor class \(K_X\) and canonical sheaf \(\om_X\), if no confusion can arise. We usually write \(K_{X}=\Oh(k)\) to represent \(K_X = kD\).   
\end{itemize}

 \subsection{Graded rings and Polarized varieties}We call a  pair \((X,D)\)   a  \emph{polarized variety}  if   \(X\) is a  normal   projective  algebraic variety and \(D\)  a \(\QQ\)-ample Weil divisor on \(X\), i.e. some integer multiple of \(D\) is a  Cartier divisor.  One gets an associated  finitely generated graded ring
\begin{displaymath}
R(X,D)=\bigoplus_{n \geq 0}H^0\left(X,\Oh_X(nD)\right).
\end{displaymath} 
It is called a \emph{projectively Gorenstein} if the ring \(R(X,D)\) is a Gorenstein ring. A surjective morphism from a free graded ring \(k[x_0, . . . , x_n]\) to $R(X,D)$    
gives the embedding 
\[i : X = \Proj R(X,D) \into \PP(a_0, \cdots , a_n)\]
where \(a_i=\deg(x_i)\) and with the divisorial sheaf \(\Oh_X(D)\) being isomorphic to \(\Oh_X(1) = i^*\Oh_{\PP}(1)\). 

 \noindent The Hilbert series of a polarized projective variety \((X,D)\)  is given by    
\begin{equation}\label{hs} P_{(X,D)}(t)=\sum_{m \geq 0}h^0(X,mD)\ t^{m},
\end{equation}where \(h^0(X,mD)=\dim H^0(X,\Oh_X(mD)).\)
We usually write \(P_X(t)\) for the Hilbert series  and by the standard 
Hilbert--Serre theorem \cite[Theorem 11.1]{atiyah}, 
$P_X(t)$ has the following  compact form
\begin{equation}
P_X(t)=\dfrac{N(t)}{ \displaystyle\prod_{i=0}^a\left(1-t^{a_i}\right)},
\label{reducedhs}
\end{equation}
where \(N(t)\) is a palindromic  polynomial of degree \(q\), as \(X\) is projectively Gorenstein.    
 \subsection{Rigid del Pezzo surfaces}
\begin{dfn}
An isolated orbifold point \(Q\) of type  \(\frac1r(a_1,\ldots,a_n)\) is the quotient of  \(\AA^n\) by the cyclic group \(\mu_r\),$$\epsilon: (x_1, \ldots,x_n)\mapsto (\epsilon^{a_1}x_1,\ldots,\epsilon^{a_n}x_n)$$
such that \(\GCD(r,a_i)=1\) for \(1\le i\le n \), \(0<a_i<r,\) and \(\epsilon\) is a primitive generator of \(\mu_r\). 
\end{dfn} 

 A \emph{del~Pezzo surface} \(X\) is a two dimensional algebraic variety with an ample anti-canonical divisor class $-K_X$. If, at worst,  $X$ contains    isolated orbifold points then we call it an  \emph{orbifold or a log del Pezzo surface}.   The \emph{Fano index} \( I\) of \(X\) is the largest positive integer $I$ such that $-K_X = I\cdot D$ for some divisor  $D $ in the divisor class group of $X$. An orbifold del Pezzo surfaces  \(X\subset \PP(a_i)  \) of codimension \(c\) is   well-formed if the singular locus of \(X\) consists of at most isolated points. It is quasismooth  if the affine cone 
\(\widetilde X =\Spec R(X,D) \subset \AA^{n+1}\)  is smooth outside its vertex \(\underline{0}\). 

A  singularity admitting a  $\QQ$-Gorenstein smoothing is called \emph{a T-singularity}\cite{Kollar-SB}.  A singularity which is rigid under  $\QQ$-Gorenstein smoothing is called  a rigid or  \emph{R-singularity}\cite{AK-SC}.  The following   
characterization of a   $T$-singularity and $R$-singularity are useful in our context \cite{dp-CP}.

\begin{dfn}
\label{dfn:T_R}
Let  $Q = \frac{1}{r}(a,b)$ be an orbifold point and take   $m = \GCD(a+b,r)$, $s = (a+b)/m$ and $k=r/m$  then $Q$ has a form  $\frac{1}{mk}(1,ms-1). $ Moreover \(Q\) is called a \emph{T}-singularity if $k \mid m$~\cite{Kollar-SB}   and an \emph{R}-singularity   if $m<k$ ~\cite{AK-SC}.
  
\end{dfn}
In the two dimensional case, any  orbifold point \(\frac1r(a,b)\) can be represented as \(\frac1r(1,a')\) by choosing  a different  primitive generator of the cyclic group \(\mu_r \) and the following Lemma follows from it.  
\begin{lema}
\label{lema-sing-relations}
Let \(Q_1=\frac1r(1,a)\) and \(Q_2=\frac1r(1,b)\) be isolated orbifold points. Then  \(Q_1=Q_2\) if and only if  \(a=b \textrm{ or }ab\equiv 1 \mod r \). 
\end{lema}
By using the    fact that each   orbifold point on a surface   can be written as \(\frac{1}{r}(1,a) \) and by applying   Lemma  \ref{lema-sing-relations}  on the  all  possible isolated rigid orbifold points of type \(\frac1r(1,a);\; 3 \le r \le 10\), we get to the  following Lemma. 
\begin{lema}\label{lema-sings}Let \(3\le r\le 10\) then any isolated rigid orbifold point \(\dfrac1r(a,b)\) is  equivalent to one of the following.
 \renewcommand*{\arraystretch}{1.5}
$$\left\{\begin{array}{ccccccc}\frac13(1,1),&\frac15(1,1),&\frac15(1,2),&\frac16(1,1),&\frac16(1,5),&\frac17(1,1),&\frac17(1,2),\\\frac17(1,3),&\frac18(1,1),&\frac18(1,5),&\frac19(1,1),&\frac19(1,4),&\frac1{10}(1,1),&\frac1{10}(1,3)\\
\end{array}\right\}$$
\end{lema}

%

 \subsection{Ambient  varieties}  In this section we briefly recall the definition of  weighted Grassmannian \(\w\Gr(2,5)\) and \(\w(\PxP\)) which we use, apart from weighted projective spaces,  as rep-quotient varieties  for the construction
of our rigid orbifold del Pezzo surfaces; following the notion introduced in \cite{dp-CH}.
\subsubsection{Weighted Grassmannian \(\w\Gr(2,5)\)} This part is wholly  based on material from   \cite[Sec.2]{wg} . 
Let $w:=(w_1,\cdots,w_5)$   be a tuple  of  all integers or  all half integers such that $$w_i+w_j>0,\; 1\le i<j\le 5,$$ Then  the quotient of the affine cone over Grassmannian minus the origin   $\widetilde{\Gr(2,5)} \backslash\{\underline 0\} $ by $\CC^\times$ given by: 
   $$\epsilon:x_{ij}\mapsto \epsilon^{w_i +w_j }x_{ij}$$ 
  is called weighted Grassmannians \(\wgr\) where \(x_{ij}\) are Pl\"ucker coordinates of the embedding \(\Gr(2,5)\into \PP\left(\bigwedge^2 \CC^5\right).\) Therefore we get the embedding
$$\wgr\into \PP\left( a_{ij}: 1\le i <j\le 5, a_{ij}=w_i+w_j\right).  $$ The image of \(\Gr(2,5)\) and \(\wgr\) under the Pl\"ucker embedding is defined by  five \(4 \times 4\) Pfaffians of  the $5\times 5$ skew symmetric matrix  
   \begin{equation}\label{eq:pf_mat}\left(\begin{matrix} 
x_{12}&x_{13}&x_{14}&x_{15}\\ 
&x_{23}&x_{24}&x_{25}\\ 
&&x_{34}&x_{35}\\ 
&&& x_{45} \end{matrix}\right)
,\end{equation}
where we only write down the upper triangular part. Explicitly the defining  equations are:
\[\Pf_i=x_{jk}x_{lm}-x_{jl}x_{km}+x_{jm}x_{lm},\]
where \(1\le j<k<l<m\le5\) are four integers and \(i\) makes up the fifth one in \(\{1,2,3,4,5\}\). 
In examples we usually write down the corresponding matrix of weights, replacing \(x_{ij} \) with \(a_{ij}\); to represent the given \(\wgr\). 
 
If $\wgr$ is wellformed then the orbifold canonical  divisor class is \begin{equation}\label{eq:K_G}K_\wgr= \left(-\frac12\sum_{1\le i<j\le 5} a_{ij} \right)D, \end{equation} for a  divisor \(D\) in the class group of \( \wgr\).
\subsubsection{Weighted \(\PxP\)} This section recalls the definition of weighted \(\PxP\) from \cite{BKQ,Sz05}. 
 Let   $b=(b_1,b_2,b_3)$ and $c=(c_1,c_2,c_3)$ be two integer or half integer vectors satisfying  $$b_1+c_1 > 0,\quad b_i\le b_j \text{ and }\quad c_i\le c_j \text{ for } 1\le i\le j \le 3,$$  and \(\Si_P\) denotes the Segre embedding $\PxP\into \PP^8(x_{ij}) $.  If \( \widetilde {\Si_P}\) is the affine of this Segre embedding, then  the weighted \(\PxP\) variety   $\w\Si_P$ is the quotient  of the punctured affine cone $\widetilde{\Si_P}\backslash\{\underline 0\}$  by  $\CC^\times$:
$$\epsilon:x_{ij}\mapsto\epsilon^{b_i+c_j}x_{ij}, \; 1\le i,j\le 3.$$
 Thus for a choice of $b,c$, written together as a single  input parameter $p=(b_1,b_2,b_3;c_1,c_2,c_3)$, we get the embedding  $$\w\Si_P\into \PP^8(a_{ij}: a_{ij}=b_i+c_j;1\le i,j\le 3). $$  The equations are defined by $2\times 2$ minors of a size 3 square matrix which we usually refer to as the weight matrix and write it as
\begin{equation}\label{eq:wtmx_P2}
\begin{pmatrix} a_{11} & a_{12} & a_{13} \\ a_{21} & a_{22} & a_{23} \\ a_{31} & a_{32} & a_{33} \end{pmatrix} \text{ where }a_{ij}=b_i+c_j; 1\le i,j\le 3.
\end{equation}

 If $\w\Si_P$ is wellformed then the  canonical divisor class is given by \begin{equation}\label{eq:K_P2}K_{\w\Si_P}=\left(-\sum_{ i=j} a_{ij} \right)D,\end{equation}
for a  divisor \(D\) in the class group of \( \w\Si_P\).\\[5mm]

\section{Proof of the main theorem}
In this section we provide details of various steps  of our calculations which together provide the proof of Theorem \ref{MainResult}. In summary,  for each codimension and   Fano index, we first search for the list of  candidate varieties using the algorithmic approach of \cite{BKZ,QJSC}. The candidate lists comes with a suggestive  basket(s) of orbifold points and invariants. Then we perform theoretical analysis of each candidate   to establish the existence or non-existence of candidate surfaces with given basket  and invariants.

\subsection{Algorithm} We briefly recall the algorithm from \cite{QJSC} which we used to compute the candidate lists of examples. The key part of it is   based on the orbifold Riemann--Roch formula
of Bukcley, Reid and Zhou \cite{BRZ} which provides a decomposition of the Hilbert series of \(X\) into a  smooth part and a singular part.
It roughly states that if \(X\) is an algebraic variety with basket \(\mathcal B=\{k_i \times Q_i: m_i \in \ZZ_{>0}\}\) of isolated orbifold points then its Hilbert series has a decomposition into a smooth part \(P_\sm(t)\) and orbifold part $\sum k_iP_{Q_i}(t)$;
\begin{equation}\label{hs-decomp}P_X(t)=P_\sm(t)+\sum k_iP_{Q_i}(t).\end{equation}
The algorithm searches for all orbifolds of fixed dimension \(n \) having  fixed orbifold canonical class \( K_{X}=\Oh(k)\) in a given ambient rep-quotient variety. Indeed, if  \(X\) is a Fano variety of index \(I\) then \(k=-I. \) The algorithm has the following steps. 
\begin{enumerate}
\item Compute the Hilbert series and orbifold canonical class of ambient rep-quotient variety.
\item 
Find all possible embeddings of \(n\)-folds \(X\) with \(\om_X=\Oh(k)\) by  applying the adjunction formula.
\item 
For each possible \(n\)-fold embedding of \(X\), compute the Hilbert series \(P_X(t)\) and  the smooth term \(P_\sm(t)\). 
\item Compute the list of all possible \(n\)-fold isolated orbifold points from the ambient weighted projective space containing \(X\).
\item For each subset of the list of possible orbifold points determine the multiplicities  \(k_{i}\) given in equation \eqref{hs-decomp} of the orbifold terms \(P_{Q_i}(t)\). 
\item If \(k_i\ge 0\) then \(X\) is a candidate \(n\)-fold with suggested basket of isolated orbifold points.

\end{enumerate}

\subsection{Bounds on search parameters}\label{Bounds}  We perform our search in the order of increasing sum of the weights on the ambient weight projective space \(\PP(a_0,\ldots,a_n)\) containing \(X\).  The search is theoretically unbounded in each codimension in two directions: there is no bound on the sum of  weights \(W=\sum a_i\)  of the ambient weighted projective space containing  \(X\) and the Fano index \(I\) is also unbounded. 

In each codimension, we at least search for  polarized rigid del Pezzo surfaces     \( X \into\PP(a_i)\) such that    \[W -I \le 50,\textrm{ for } 1 \le I \le10.\] If the last candidate example is found for the adjunction number   \(q=W-I  \) of the Hilbert numerator \(N(t),  \) then we further search for all possible cases such that \[W -I \le N \textrm{ where } N= \Max(2q,50),\] to absolutely  minimize the possibility of any missing examples.
 Similarly  in each codimension if we find the last example in search domain  \( W -I \le 50\) for index \(I > 5\) then we search for examples  up to index  \(2I\).  For example,  in the hypersurface case the maximum
value of   \(I\) across all candidates was 8, so we searched until index 16 in this case. Similarly, for index 2 hypersurfaces  we got the last candidate when \( W -2 = 36\) so we searched for all cases with \(W -2 \le 72.\) 
\begin{table}[h]

\caption{ The following table summarises the number of surfaces we obtained for each Fano index \(I \) in each codimension and exact search domain in each case.  First column contains the codimension of each surface and the rest of the columns contain a pair of numbers. First number is the number of examples of given index and the second one gives the maximum value of \(q=W_{\max}-I\) for which the last candidate surface was found; the classification is  complete till \( N=\Max(50,2q)\). The entries with no second number means that no examples were found for \(q\le 50.\)    
\label{tab-summary}}
\begin{tabular}{ |l|l|l|l|l|l|l|l|l |l|}
  \hline
 \multirow{2}{*}{Codimension} &\multicolumn{9}{|c|}{Fano Index, (\(q \))} \\\cline{2-10}
  &1&2&3&4&5&6&7&8&9--16 \\
   \hline
\evnrow   1 &11(28)&44(36) &6(15)&6(21)&6(21)&2(16)&2(17)&4(15)&0\\ \hline
   2 &15(22)&8 (29)&1(26)&1(22)&0&0&0&0&\\ \hline
\evnrow   3 &12(33)&7 (43)&1(19)&1(26)&0&0&0&0&\\ \hline
   4 &12(42)&6 (48) &0&0&1(30)&0&0&1(42)&0\\ \hline
\end{tabular}
\end{table}

 \subsection{Computing invariants}
\label{Comp-Invariants}  We describe how we calculate each of the following \(qG\)-deformation invariants appearing in  tables \ref{tab-hypersurface}--\ref{tab-P2xP2}. 
\begin{enumerate}
\item  \textbf{First plurigenus} \(h^0(-K_X)\):  If it is equal to zero then we can easily conclude that   \(X\) does not admit a \(qG\)-deformation to a toric variety and such surfaces are not of class \(TG\). We compute it as  the coefficient of \(t^{I}\) in the Hilbert series \eqref{hs} where \(I\) is the Fano index of \(X\). 
\item \textbf{Intersection number} \(-K_X^2\):   It can be defined as an anti-canonical degree of \(X\) which we  calculate from the Hilbert series \(P_X(t)\) of \(X\). In a surface case  
\[P_X(t)=\dfrac{H(t)}{(1-t)^3},\] where \(H(t)  \) is a rational function with only positive coefficients. Then  for a generic divisor \(D\) in the class group, we have   \(D^2=H(1)\). Consequently for an  orbifold del Pezzo surface of index \(I\), we have \(-K_X^2= I^2D^2\). 
 \item \textbf{Euler Characteristics} \(e(X)\): We were  able to compute the Euler characteristics of \(X\)  in   hypersurface and  complete intersection cases by using  Blache's formula \cite[2.11-14]{Blache};
\begin{equation}
e(X)=e_\orb(X)+\displaystyle\sum_{r(Q)\in\mathcal B} \frac{r-1}{r}
\end{equation}
where \(r\) is the local index of each orbifold point. It was
 applied in the Appendix of \cite{BF} to illustrate the computation for a  hypersurface. The formula has natural generalization to the cases of complete intersections \[X_{d_1,\ldots,d_k}\subset \PP(a_0,\ldots,a_n)\]in higher codimension. We can computer \(e_\orb(X)\) as:
\begin{equation} 
e_\orb(X)=\text{coefficient of } t^{n-k} \text{ in the series expansion of}\left(\dfrac{\prod_{}(1+a_it)}{\prod (1+d_it)}\deg(X)\right). 
\end{equation}
\item  \textbf{Picard rank  \(\rho(x) \)}:  We were able to calculate it explicitly when \(X\) is a hypersurface in \( \PP^{3}(a_{i})\) by using \cite[4.4.1]{dolgachev}. Given a  hypersurface \[X_{d }\into \PP(a_0,a_1,a_2,a_3),\]  let  
\[l=\text{coefficient of } t^{2d-\sum a_i}\ \text{ in the series  expansion of}\left(\prod\dfrac{t^{d-a_i}-1}{ t^{a_i}-1}\right),\]  then \(\rho(X)=l+1.\) In cases of complete intersection examples we were able to identify those examples which are not prime, i.e. the Picard rank greater than 1.
From \cite{dp-Dong}, we know that  if the Picard rank of a log del Pezzo surface is 1 then \(0< e_\orb(X)\le 3\).  Therefore  for each codimension  2 complete intersection in Table \ref{tab-CI},  we list  \(e_\orb(X)\) and those with \(e_\orb(X) >3\)  have Picard rank greater than 1.
\end{enumerate}
\subsection{Theoretical singularity  analysis}
\label{sing-analysis}

 The last step of the calculation is the theoretical singularity analysis of each candidate orbifold. We prove  that the general member  \(X\) in each family is wellformed and quasismooth. We first compute the dimensions of intersection of all orbifold strata with \(X\) to establish that \(X\) is wellformed. This should be less than or equal to zero for a surface to be wellformed, i.e. it does not contain any singular lines. \ 

 The next  step is  to show that \(X\) is quasismooth. It is not so difficult when \(X\) is a hypersurface or complete intersection: one can use the criteria given in \cite[ 8]{fletcher}. In cases of  codimension 3 and 4 examples, we consider \(X\) as complete intersections  in $\wgr$  or  in the Segre embedding of weighted $\PxP$  or in some projective cone(s) over either of those ambient  varieties. So \(X\) may not only have singularities from the ambient weighted projective but  it   may also contain singularities on the base loci of linear systems of the intersecting weighted homogeneous forms. In such cases we mostly  prove the quasismoothness on the base locus by using computer algebra system \textsc{magma} \cite{magma}. We write down  explicit equations for \(X\) over the rational numbers and show that it is smooth, see  \cite[2.3]{QMOC} for more details. To prove quasismoothness on an orbifold point \(Q\) of type \(\frac1r(a,b)\), which is mostly a coordinate point corresponding to some variables \(x_i\) with \(\deg(x_i)=r\), we proceed as follows. If $c$ is the codimension of $X$ then we find $c$ tangent variables  $x_m$ \cite{BZ}, i.e.  we find c polynomials having a monomial of type
$x_i^lx_m$. We can locally remove these variables by  using the implicit function theorem. Moreover, if  two other  variables  have  weights  $a$ and $b$ modulo $r$  then \(Q\) is a quasismooth point of type $ \frac1r (a, b)$.

\section{Sample Calculations}
In this section we provide sample calculations of examples  given in  tables \ref{tab-hypersurface}--\ref{tab-P2xP2}.  
\begin{eg}\(\#81\) Consider the weighted projective space \(\PP(1,5,7,10)\)  with variables \(x,y,z \text{ and } w\) respectively, then the canonical class \(K_\PP=\Oh(-23)\). The  generic weighted homogenous polynomial of degree 15, \[f_{15}=k_1~x^{15}+k_2~y^3+k_3~yw+k_4~xz^2+\cdots, \quad k_i\in \CC;\] defines a del Pezzo surface \(X_{15}\into \PP(x,y,z,w)\) of Fano index 8, i.e. \(K_X=\Oh(-8)\). The polynomial \(f_{15}\) does not contain monomials of pure power in \(w\) and \(z\) so \(X\) contains the orbifold points \(p_{1}=(0,0,0,1)\) and \(p_{2}=(0,0,1,0)\). By applying the implicit function theorem we can remove the variable \(y\) near  the point \(p_{1}\) by using the monomial \(yw\) and \(x,z\) are local variables near this point. Therefore \(X\) contains an orbifold point of type \(\frac1{10}(1,7)=\frac{1}{10}(1,3)\)(Lemma \ref{lema-sing-relations}).   Similarly near \(p_2\) the local variables are \(y\) and \(w\), so we get an orbifold point of type \[\frac17(5,10)=\frac17(3,5)=\frac17(1,4)=\frac17(1,2).\]
The coordinate point of weight 5 does not lie on \(X\) but one dimensional singular stratum \(\PP^1(y,w)\) intersects with \(X\) non-trivially and by \cite[Lemma~9.4]{fletcher} the intersection is in  two points. One of them is \(p_1\) and the other   can be taken as \(p_3=(0,1,0,0)\) which corresponds to weight 5 variable. By using the above arguments we can show that it is a singular point of type \(\frac15(1,2)\). Thus \(X\) contains exactly the  same basket of singularities as given by the computer search and it is a wellformed and quasismooth rigid del Pezzo surface of Fano index \(8\). Moreover the vector space \[H^{0}(X,-K_X)=H^0(X,8D)=<x^8,x^3y,xz>,\]
so \(h^0(-K_X)=3\).
\end{eg}
\begin{eg}\(\#126\) Consider the weighted Grassmannian \(\wgr\) 
\[\wgr\into \PP(1^2,3^3,5^4,7) \text{ with weight matrix   } \begin{pmatrix}1&1&3&3\\&3&5&5\\&&5&5\\&&&7\end{pmatrix},\]
Then by equation \eqref{eq:K_G} the canonical divisor class \(K_{\wgr}=\Oh(-19)\). The weighted complete intersection of \( \wgr \) with two forms of degree 3 and two forms of degree 5;  \[X=\wgr \cap(f_3)\cap(g_3)\cap(f_5)\cap (g_5)\into \PP_{(x_1,x_2,y_1,z_1,z_2,w_1)}(1^2,3,5^2,7) \] is a del Pezzo surface with \(K_X=\Oh(-19+(3+3+5+5))=\Oh(-3)\).  We can take  \(X\) to be defined by the maximal Pfaffians of 
 \begin{equation}\label{eq:ex_pf}\begin{pmatrix}x_1&x_2&f_3&g_3\\&y_{1}&f_{5}&g_{5}\\&&z_1&z_2\\&&&w_1\end{pmatrix},\end{equation}
 where \(f_3,g_3,f_5 \text{ and } g_5\) are general weighted homogeneous forms in given variables and they remove the variables of the corresponding degrees from the ambient \( \wgr\). The coordinate point corresponding to  \(w_1\) lies on \(X\). From the equations we have \(x_1,x_2 \text{ and } y_1\) as tangent variables and \(z_1,z_2\) as local variables. Therefore it is an orbifold point of type \(\frac17(5,5)=\frac17(1,1). \) The locus \(X \cap \PP(5,5)\) is locally  a quadric in \(\PP^1\) which defines two points.  By similar application of implicit function theorem we can show that each  is an orbifold point of type \(\frac15(1,2).\)   The restriction of \(X\) to weight \(3\) locus is an empty set, so \(X\) contains no further orbifold points.  
 To show the quasismoothness on the base locus we use the computer algebra and write down equations for \(X\). For example, if we choose 
 \renewcommand*{\arraystretch}{1.5}
\[\begin{array}{cc}f_3=3x_1^3 + 3x_2^3,  &f_5= x_2^5 + x_1^2y_1 + x_2^2y_1 + z_1 + z_2,\\g_3= x_2^3 + y_1,  &g_5=x_1^5  + 2x_1^2y_1 + 3x_2^2y_1  + 3z_2 \end{array}\]
then the Pfaffians of \eqref{eq:ex_pf} gives a quasismooth surface. Thus \(X\) is an orbifold del Pezzo surface of Fano index 3 with singular points; \(2\times \frac15(1,2)\) and \(\frac17(1,1)\).    
\end{eg}
As we mentioned in  section \ref{sing-analysis}  that we prove the existence of given orbifold del Pezzo surface by theoretical singularity analysis. Then only those which are quasismooth, wellformed and  having correct basket of singularities appear in tables of examples. There are total 8 candidate examples which fails to be quasismooth and we discuss one of them below in detail. No candidate example fails for not being  wellformed.  

\begin{eg}(Non working candidate)
A computer search also gives a candidate complete intersection orbifold del Pezzo surface of Fano index 2 given by\[X_{6,30}\into\PP_{(x,y,z,t,u)}(1,3,9,10,15).\]
Then \(F_6=f(x,y)\) (since other variables have weight higher than 6) and \[F_{30}=x^{30}+x^{27}y+yz^3+\cdots \] are the defining equations of \(X. \) The coordinate point $p=(0,0,1,0,0)$ lies on \(X\) as no pure power of \(z\) appear in \(F_{30}\). Now we can not find two tangent variables to \(z\) in the equations of \(X \) which implies that the rank of the Jacobian matrix of \(X\) at \(p\) is equal to 1 which is  less than its codimension, so \(X \) is not quasismooth at \(p.\) Thus \(X\) is a del Pezzo surface which is not quasismooth and does not appear in the following tables. \end{eg}

\noindent\textbf{Concluding remark}: One can use this approach to construct and classify orbifold del Pezzo surfaces with any quotient singularity in a given fixed format, under certain bounds. Moreover, we can also construct examples with rigid orbifold points of type \(\frac 1r (1,a)\) for \(r\ge 11\) but as the weights higher the computer search output becomes slower due to the nature of  algorithm. Therefore, we restrict ourself to the cases with \(r\le 10\).   
\section*{Acknowledgements} I would like to thank Erik Paemurru for  pointing me to \cite{dolgachev} for computations of Picard rank in the hypersurface case. I am grateful to an anonymous referee for their feedback which improved the exposition of the paper significantly. I would  like to acknowledge the support provided by Deanship of Scientific Research (DSR) at King Fahd University of
Petroleum and Minerals for funding this work through project No. SB191029.  
\bibliographystyle{amsplain}
\bibliography{References}

\providecommand{\bysame}{\leavevmode\hbox to3em{\hrulefill}\thinspace}
\providecommand{\MR}{\relax\ifhmode\unskip\space\fi MR }
\providecommand{\MRhref}[2]{%
  \href{http://www.ams.org/mathscinet-getitem?mr=#1}{#2}
}
\providecommand{\href}[2]{#2}
\begin{thebibliography}{10}

\bibitem{dp-AMS}
Mohammad Akhtar, Tom Coates, Alessio Corti, Liana Heuberger, Alexander
  Kasprzyk, Alessandro Oneto, Andrea Petracci, Thomas Prince, and Ketil
  Tveiten, \emph{Mirror symmetry and the classification of orbifold {D}el
  {P}ezzo surfaces}, Proceedings of the American Mathematical Society
  \textbf{144} (2016), no.~2, 513--527.

\bibitem{AK-SC}
Mohammad Akhtar and Alexander Kasprzyk, \emph{Singularity content}, arXiv
  preprint arXiv:1401.5458 (2014).

\bibitem{atiyah}
M.~F. Atiyah and I.~G. Macdonald, \emph{Introduction to commutative algebra},
  Addison-Wesley Publishing Co., Reading, Mass.-London-Don Mills, Ont., 1969.

\bibitem{Blache}
Raimund Blache, \emph{Chern classes and {H}irzebruch--{R}iemann--{R}och theorem
  for coherent sheaves on complex-projective orbifolds with isolated
  singularities}, Mathematische Zeitschrift \textbf{222} (1996), no.~1, 7--57.

\bibitem{magma}
Wieb Bosma, John Cannon, and Catherine Playoust, \emph{The {M}agma algebra
  system. {I}. {T}he user language}, J. Symbolic Comput. \textbf{24} (1997),
  no.~3-4, 235--265, Computational algebra and number theory (London, 1993).

\bibitem{BKQ}
G.~Brown, A.~M. Kasprzyk, and M.~I. Qureshi, \emph{Fano 3-folds in
  $\mathbb{P}^2\times \mathbb{P}^2$ format, {Tom} and {J}erry}, European
  Journal of Mathematics \textbf{4} (2018), no.~1, 51--72.

\bibitem{BF}
Gavin Brown and Enrico Fatighenti, \emph{Hodge numbers and deformations of
  {F}ano 3-folds}, arXiv preprint arXiv:1707.00653 (2017).

\bibitem{BKZ}
Gavin Brown, Alexander~M Kasprzyk, and Lei Zhu, \emph{Gorenstein formats,
  canonical and {C}alabi--{Y}au threefolds}, Experimental Mathematics (2019),
  1--19.

\bibitem{BZ}
Gavin Brown and Francesco Zucconi, \emph{Graded rings of rank 2 {S}arkisov
  links}, Nagoya Mathematical Journal \textbf{197} (2010), 1--44.

\bibitem{BRZ}
A.~Buckley, M.~Reid, and S.~Zhou, \emph{Ice cream and orbifold
  {R}iemann--{R}och}, Izvestiya: Mathematics \textbf{77:3} (2013), 461--486.

\bibitem{dp-CP}
D~Cavey and T~Prince, \emph{Del {P}ezzo surfaces with a single $\frac{1}{k} (1,
  1)$ singularity}, Journal of the Mathematical Society of Japan (2020).

\bibitem{Fano-CK}
Daniel Cavey and Edwin Kutas, \emph{Classification of minimal polygons with
  specified singularity content}, arXiv preprint arXiv:1703.05266 (2017).

\bibitem{dp-exceptional}
Ivan Cheltsov, Jihun Park, and Constantin Shramov, \emph{Exceptional {D}el
  {P}ezzo hypersurfaces}, Journal of Geometric Analysis \textbf{20} (2010),
  no.~4, 787--816.

\bibitem{dp-zoo}
Ivan Cheltsov and Constantin Shramov, \emph{{D}el {P}ezzo zoo}, Experimental
  Mathematics \textbf{22} (2013), no.~3, 313--326.

\bibitem{CCGGK}
Tom Coates, Alessio Corti, Sergey Galkin, Vasily Golyshev, and Alexander
  Kasprzyk, \emph{Mirror symmetry and {F}ano manifolds}, Proceedings of
  European Congress of Mathematics (Krakow, 2-7 July, 2012), January 2014 (824
  pages), ISBN 978-3-03719-120, DOI 10.4171/120-1/16 (2012).

\bibitem{wg}
A.~Corti and M.~Reid, \emph{Weighted {G}rassmannians}, Algebraic geometry
  (M.~C. Beltrametti, F.~Catanese, C.~Ciliberto, A.~Lanteri, and C.~Pedrini,
  eds.), de Gruyter, Berlin, 2002, pp.~141--163.

\bibitem{dp-CH}
Alessio Corti and Liana Heuberger, \emph{{D}el {P}ezzo surfaces with
  {$\frac{1}{3}(1,1)$} points}, Manuscripta Math. \textbf{153} (2017), no.~1-2,
  71--118.

\bibitem{dolgachev}
I.~Dolgachev, \emph{Weighted projective spaces}, Group actions and vector
  fields, Lec Note in Mathematics, Springer-Verlag,956, 1981, pp.~34--71.

\bibitem{dp-Dong}
DongSeon Hwang, \emph{On the orbifold euler characteristic of log del {P}ezzo
  surfaces of rank one}, Journal of the Korean Mathematical Society \textbf{51}
  (2014), no.~4, 867--879.

\bibitem{fletcher}
A.~R. Iano-Fletcher, \emph{Working with weighted complete intersections},
  Explicit Birational Geometry of 3-folds, vol. 281, London Math. Soc. Lecture
  Note Ser, CUP, 2000, pp.~101--173.

\bibitem{dp-KNP}
Alexander Kasprzyk, Benjamin Nill, and Thomas Prince, \emph{Minimality and
  mutation-equivalence of polygons}, Forum of Mathematics, Sigma, vol.~5,
  Cambridge University Press, 2017.

\bibitem{Kollar-SB}
J{\'a}nos Koll{\'a}r and Nicholas~I Shepherd-Barron, \emph{Threefolds and
  deformations of surface singularities}, Inventiones mathematicae \textbf{91}
  (1988), no.~2, 299--338.

\bibitem{Mayanskiy}
Evgeny Mayanskiy, \emph{Weighted complete intersection del {P}ezzo surfaces},
  arXiv preprint arXiv:1608.02049 (2016).

\bibitem{dp-Miura}
Takayuki Miura, \emph{Classification of del {P}ezzo surfaces with
  $\frac{1}{3}(1,1)$ and $\frac{1}{4}(1, 1) $ singularities}, arXiv preprint
  arXiv:1903.00679 (2019).

\bibitem{dp-Erik}
Erik Paemurru, \emph{Del {P}ezzo surfaces in weighted projective spaces},
  Proceedings of the Edinburgh Mathematical Society \textbf{61} (2018), no.~2,
  545--572.

\bibitem{QS}
M.~I. Qureshi and B.~Szendr{\H o}i, \emph{Constructing projective varieties in
  weighted flag varieties}, Bull. Lon. Math Soc. \textbf{43} (2011), no.~2,
  786--798.

\bibitem{QS-AHEP}
M.~I. Qureshi and Bal{\'a}zs Szendr{\H{o}}i, \emph{Calabi-{Y}au threefolds in
  weighted flag varieties}, Adv. High Energy Phys. (2012), Art. ID 547317, 14
  pp.

\bibitem{QJSC}
Muhammad~Imran Qureshi, \emph{Computing isolated orbifolds in weighted flag
  varieties}, Journal of Symbolic Computation \textbf{79, Part 2} (2017), 457
  -- 474.

\bibitem{QMOC}
Muhammad~Imran Qureshi, \emph{Biregular models of log del {P}ezzo surfaces with
  rigid singularities}, Mathematics of Computation \textbf{88} (2019), no.~319,
  2497--2521.

\bibitem{dp-RS}
Miles Reid and Kaori Suzuki, \emph{Cascades of projections from log del {P}ezzo
  surfaces}, London Mathematical Society Lecture Note Series, p.~227–250,
  Cambridge University Press, 2004.

\bibitem{Sz05}
Bal{\'a}zs Szendr{\H o}i, \emph{On weighted homogeneous varieties}, Unpublished
  manuscript, 2005.

\end{thebibliography}

\appendix
\section{Table of Examples}

 \subsection*{Notations in Tables}
 \begin{itemize}
\item The column \(X\) represents a del Pezzo surface and the corresponding  weighted projective space containing
\(X\);  the subscripts give the  equation degrees of \(X\). The column \(I\) lists the Fano index of \(X\).
\item  The next two columns contain the anti-canonical degree \(-K_X^2\) and the first plurigenus \(h^0(-K_X)\). If \(h^0(-K_X)=0\) the \(X\) is not of class \(TG.\)
\item \(e(X)\) denotes the topological Euler characteristics of \(X,~\rho(X)\) is the rank of Picard group of \(X\), and \(e_{\orb}(X)\) denotes the orbifold Euler number of \(X\). \(\rho(X)\) is only listed in table \ref{tab-hypersurface} of hypersurfaces and \(e_{\orb}(X)\) only in table \ref{tab-CI} of complete intersections, as discussed in \ref{Comp-Invariants}. 
\item The column \(\sB\) represents the basket of singular points of \(X\).
\item In Table \ref{tab-Pf} and \ref{tab-P2xP2}, the last column represents the matrix of weights, which provides weights of ambient weighted projective space containing \(\wgr\) or weighted \(\PxP\) variety.
\item We provide references to those examples which appeared in  \cite{dp-CH} and \cite{Fano-CK}, primarily in a toric setting. 
\end{itemize}  

 \renewcommand*{\arraystretch}{1.5}
\begin{longtable}{>{\hspace{0.5em}}llcccccc<{\hspace{0.5em}}}
\caption{ Hypersurfaces in \(w\PP^3\) }
\label{tab-hypersurface}\\
\toprule
\multicolumn{1}{c}{S.No}&\multicolumn{1}{c}{$X$ }&\multicolumn{1}{c}{$I$}&\multicolumn{1}{c}{$-K_X^2$}&\multicolumn{1}{c}{$h^0(-K_X)$}&\multicolumn{1}{c}{$e(X)$}&\multicolumn{1}{c}{${\rho}(X)$}&\multicolumn{1}{c}{Basket $\mathcal{B}$}\\
\cmidrule(lr){1-1}\cmidrule(lr){2-2}\cmidrule(ll){3-5}\cmidrule(lr){6-7}\cmidrule(lr){8-8}
\endfirsthead
\multicolumn{7}{l}{\vspace{-0.25em}\scriptsize\emph{\tablename\ \thetable{} continued from previous page}}\\
\midrule
\endhead
\multicolumn{7}{r}{\scriptsize\emph{Continued on next page}}\\
\endfoot
\bottomrule
\endlastfoot
\oddrow $\rownumber$ &$X_{3}\subset \PP(1^4) $& 1 & $3$ & $4$& $9$ & $7$ & $ { }
$ \\

  \evnrow $\rownumber$ &$X_{4}\subset \PP(1^3, 2) $& 1 & $2$ & $3$& $10$ & $8$ &
$ { } $ \\

  \oddrow $\rownumber$ &$X_{6}\subset \PP(1^2, 2, 3) $& 1 & $1$ & $2$& $11$ & $9$
& $ { } $ \\

  \evnrow $\rownumber$ &$X_{10}\subset \PP(1, 2, 3, 5) $& 1 & $1/3$ & $1$ & $11$
& $9$ & $ { \frac{1}{3}(1, 1)} $~~\cite{dp-CH} \\

 \oddrow $\rownumber$ &$X_{12}\subset \PP(2, 3^2, 5) $& 1 & $2/15$ & $0$ & $10$
& $8$ & $ { 4 \times  \frac{1}{3}(1, 1), \frac{1}{5}(1, 1) } $ \\

 \evnrow $\rownumber$ &$X_{15}\subset \PP(1, 3, 5, 7) $& 1 & $1/7$ & $1$ & $11$
& $9$ & $ { \frac{1}{7}(1, 2) } $ \\

 \oddrow $\rownumber$ &$X_{15}\subset \PP(3^2, 5^2) $& 1 & $1/15$ & $0$ & $11$ &
$9$ & $ { 5 \times  \frac{1}{3}(1, 1), 3 \times  \frac{1}{5}(1, 1) } $ \\
\evnrow $\rownumber$ &$X_{16}\subset \PP(1, 3, 5, 8) $& 1 & $2/15$ & $1$ & $12$
& $10$ & $ { \frac{1}{3}(1, 1), \frac{1}{5}(1, 1) } $ \\

 \oddrow $\rownumber$ &$X_{18}\subset \PP(2, 3, 5, 9) $& 1 & $1/15$ & $0$ & $9$
& $7$ & $ { 2 \times  \frac{1}{3}(1, 1), \frac{1}{5}(1, 2) } $ \\

 \evnrow $\rownumber$ &$X_{20}\subset \PP(2, 5^2, 9) $& 1 & $2/45$ & $0$ & $10$
& $8$ & $ { 4 \times  \frac{1}{5}(1, 2), \frac{1}{9}(1, 1) } $ \\

 \oddrow $\rownumber$ &$X_{28}\subset \PP(3, 5, 7, 14) $& 1 & $2/105$ & $0$ &
$8$ & $6$ & $ { \frac{1}{3}(1, 1), \frac{1}{5}(1, 2), 2 \times  \frac{1}{7}(1, 2) } $ \\
\hline\hline
 \evnrow $\rownumber$ &$X_{2}\subset \PP(1^4) $& 2 & $8$ & $9$& $4$ & $2$ & $ { }
$ \\

  \oddrow $\rownumber$ &$X_{4}\subset \PP(1^3, 3) $& 2 & $16/3$ & $6$ & $6$ &
$4$ & $ { \frac{1}{3}(1, 1) } $~\cite{dp-CH} \\

 \evnrow $\rownumber$ &$X_{6}\subset \PP(1^3, 5) $& 2 & $24/5$ & $6$ & $8$ & $6$
& $ { \frac{1}{5}(1, 1) } $  \cite{Fano-CK}\\

 \oddrow $\rownumber$ &$X_{6}\subset \PP(1^2, 3^2) $& 2 & $8/3$ & $3$ & $8$ &
$6$ & $ { 2 \times  \frac{1}{3}(1, 1) } $~\cite{dp-CH} \\

 \evnrow $\rownumber$ &$X_{7}\subset \PP(1^3, 6) $& 2 & $14/3$ & $6$ & $9$ & $7$
& $ { \frac{1}{6}(1, 1)  } $\cite{Fano-CK} \\

 \oddrow $\rownumber$ &$X_{8}\subset \PP(1^2, 3, 5) $& 2 & $32/15$ & $3$ & $10$
& $8$ & $ { \frac{1}{3}(1, 1), \frac{1}{5}(1, 1) } $ \\

 \evnrow $\rownumber$ &$X_{8}\subset \PP(1^3, 7) $& 2 & $32/7$ & $6$ & $10$ &
$8$ & $ { \frac{1}{7}(1, 1) } $ \\

 \oddrow $\rownumber$ &$X_{9}\subset \PP(1^2, 3, 6) $& 2 & $2$ & $3$ & $11$ &
$9$ & $ { \frac{1}{3}(1, 1), \frac{1}{6}(1, 1) } $ \cite{Fano-CK}\\

 \evnrow $\rownumber$ &$X_{9}\subset \PP(1^3, 8) $& 2 & $9/2$ & $6$ & $11$ & $9$
& $ { \frac{1}{8}(1, 1) } $ \\

 \oddrow $\rownumber$ &$X_{10}\subset \PP(1^2, 3, 7) $& 2 & $40/21$ & $3$ & $12$
& $10$ & $ { \frac{1}{3}(1, 1), \frac{1}{7}(1, 1) } $ \\

 \evnrow $\rownumber$ &$X_{10}\subset \PP(1^3, 9) $& 2 & $40/9$ & $6$ & $12$ &
$10$ & $ { \frac{1}{9}(1, 1) } $ \\

 \oddrow $\rownumber$ &$X_{10}\subset \PP(1^2, 5^2) $& 2 & $8/5$ & $3$ & $12$ &
$10$ & $ { 2 \times  \frac{1}{5}(1, 1) } $ \cite{Fano-CK} \\

 \evnrow $\rownumber$ &$X_{11}\subset \PP(1^3, 10) $& 2 & $22/5$ & $6$ & $13$ &
$11$ & $ { \frac{1}{10}(1, 1) } $ \\

 \oddrow $\rownumber$ &$X_{11}\subset \PP(1^2, 5, 6) $& 2 & $22/15$ & $3$ & $13$
& $11$ & $ { \frac{1}{5}(1, 1), \frac{1}{6}(1, 1) } $ \\

 \evnrow $\rownumber$ &$X_{11}\subset \PP(1^2, 3, 8) $& 2 & $11/6$ & $3$ & $13$
& $11$ & $ { \frac{1}{3}(1, 1), \frac{1}{8}(1, 1) } $ \\

 \oddrow $\rownumber$ &$X_{12}\subset \PP(1^2, 5, 7) $& 2 & $48/35$ & $3$ & $14$
& $12$ & $ { \frac{1}{5}(1, 1), \frac{1}{7}(1, 1) } $ \\

 \evnrow $\rownumber$ &$X_{12}\subset \PP(1^2, 6^2) $& 2 & $4/3$ & $3$ & $14$ &
$12$ & $ { 2 \times  \frac{1}{6}(1, 1) } $\cite{Fano-CK} \\

 \oddrow $\rownumber$ &$X_{12}\subset \PP(1^2, 3, 9) $& 2 & $16/9$ & $3$ & $14$
& $12$ & $ { \frac{1}{3}(1, 1), \frac{1}{9}(1, 1) } $ \\

 \evnrow $\rownumber$ &$X_{13}\subset \PP(1^2, 5, 8) $& 2 & $13/10$ & $3$ & $15$
& $13$ & $ { \frac{1}{5}(1, 1), \frac{1}{8}(1, 1) } $ \\

 \oddrow $\rownumber$ &$X_{13}\subset \PP(1^2, 3, 10) $& 2 & $26/15$ & $3$ &
$15$ & $13$ & $ { \frac{1}{3}(1, 1), \frac{1}{10}(1, 1) } $ \\

 \evnrow $\rownumber$ &$X_{13}\subset \PP(1^2, 6, 7) $& 2 & $26/21$ & $3$ & $15$
& $13$ & $ { \frac{1}{6}(1, 1), \frac{1}{7}(1, 1) } $ \\

 \oddrow $\rownumber$ &$X_{14}\subset \PP(1^2, 6, 8) $& 2 & $7/6$ & $3$ & $16$ &
$14$ & $ { \frac{1}{6}(1, 1), \frac{1}{8}(1, 1) } $ \\

 \evnrow $\rownumber$ &$X_{14}\subset \PP(1^2, 7^2) $& 2 & $8/7$ & $3$ & $16$ &
$14$ & $ { 2 \times  \frac{1}{7}(1, 1) } $ \\

 \oddrow $\rownumber$ &$X_{14}\subset \PP(1^2, 5, 9) $& 2 & $56/45$ & $3$ & $16$
& $14$ & $ { \frac{1}{5}(1, 1), \frac{1}{9}(1, 1) } $ \\

 \evnrow $\rownumber$ &$X_{15}\subset \PP(1^2, 6, 9) $& 2 & $10/9$ & $3$ & $17$
& $15$ & $ { \frac{1}{6}(1, 1), \frac{1}{9}(1, 1) } $ \\

 \oddrow $\rownumber$ &$X_{15}\subset \PP(1^2, 5, 10) $& 2 & $6/5$ & $3$ & $17$
& $15$ & $ { \frac{1}{5}(1, 1), \frac{1}{10}(1, 1) } $ \\

 \evnrow $\rownumber$ &$X_{15}\subset \PP(1, 3, 6, 7) $& 2 & $10/21$ & $1$ &
$11$ & $9$ & $ { 2 \times  \frac{1}{3}(1, 1), \frac{1}{6}(1, 1), \frac{1}{7}(1, 2) } $ \\

 \oddrow $\rownumber$ &$X_{15}\subset \PP(1^2, 7, 8) $& 2 & $15/14$ & $3$ & $17$
& $15$ & $ { \frac{1}{7}(1, 1), \frac{1}{8}(1, 1) } $ \\

 \evnrow $\rownumber$ &$X_{16}\subset \PP(1^2, 6, 10) $& 2 & $16/15$ & $3$ &
$18$ & $16$ & $ { \frac{1}{6}(1, 1), \frac{1}{10}(1, 1) } $ \\

 \oddrow $\rownumber$ &$X_{16}\subset \PP(1^2, 7, 9) $& 2 & $64/63$ & $3$ & $18$
& $16$ & $ { \frac{1}{7}(1, 1), \frac{1}{9}(1, 1) } $ \\

 \evnrow $\rownumber$ &$X_{16}\subset \PP(1^2, 8^2) $& 2 & $1$ & $3$ & $18$ &
$16$ & $ { 2 \times  \frac{1}{8}(1, 1) } $ \\

 \oddrow $\rownumber$ &$X_{17}\subset \PP(1, 3, 7, 8) $& 2 & $17/42$ & $1$ &
$11$ & $9$ & $ { \frac{1}{3}(1, 1), \frac{1}{7}(1, 1), \frac{1}{8}(1, 5) } $ \\

 \evnrow $\rownumber$ &$X_{17}\subset \PP(1^2, 7, 10) $& 2 & $34/35$ & $3$ &
$19$ & $17$ & $ { \frac{1}{7}(1, 1), \frac{1}{10}(1, 1) } $ \\

 \oddrow $\rownumber$ &$X_{17}\subset \PP(1^2, 8, 9) $& 2 & $17/18$ & $3$ & $19$
& $17$ & $ { \frac{1}{8}(1, 1), \frac{1}{9}(1, 1) } $ \\

 \evnrow $\rownumber$ &$X_{18}\subset \PP(1^2, 8, 10) $& 2 & $9/10$ & $3$ & $20$
& $18$ & $ { \frac{1}{8}(1, 1), \frac{1}{10}(1, 1) } $ \\

 \oddrow $\rownumber$ &$X_{18}\subset \PP(1^2, 9^2) $& 2 & $8/9$ & $3$ & $20$ &
$18$ & $ { 2 \times  \frac{1}{9}(1, 1) } $ \\

 \evnrow $\rownumber$ &$X_{19}\subset \PP(1^2, 9, 10) $& 2 & $38/45$ & $3$ &
$21$ & $19$ & $ { \frac{1}{9}(1, 1), \frac{1}{10}(1, 1) } $ \\

 \oddrow $\rownumber$ &$X_{20}\subset \PP(1^2, 10^2) $& 2 & $4/5$ & $3$ & $22$ &
$20$ & $ { 2 \times  \frac{1}{10}(1, 1) } $ \\

 \evnrow $\rownumber$ &$X_{21}\subset \PP(3, 6, 7^2) $& 2 & $2/21$ & $0$ & $9$ &
$7$ & $ { 3 \times  \frac{1}{3}(1, 1), \frac{1}{6}(1, 1), 3 \times  \frac{1}{7}(1, 2) } $ \\

 \oddrow $\rownumber$ &$X_{21}\subset \PP(1, 3, 9, 10) $& 2 & $14/45$ & $1$ &
$13$ & $11$ & $ { 2 \times  \frac{1}{3}(1, 1), \frac{1}{9}(1, 1), \frac{1}{10}(1, 3)) } $ \\

 \evnrow $\rownumber$ &$X_{22}\subset \PP(1, 5, 7, 11) $& 2 & $8/35$ & $1$ &
$10$ & $8$ & $ { \frac{1}{5}(1, 1), \frac{1}{7}(1, 3) } $ \\

 \oddrow $\rownumber$ &$X_{24}\subset \PP(3, 7, 8^2) $& 2 & $1/14$ & $0$ & $7$ &
$5$ & $ { \frac{1}{7}(1, 1), 3 \times  \frac{1}{8}(1, 5) } $ \\

 \evnrow $\rownumber$ &$X_{30}\subset \PP(3, 9, 10^2) $& 2 & $2/45$ & $0$ & $9$
& $7$ & $ { 3 \times  \frac{1}{3}(1, 1), \frac{1}{9}(1, 1), 3 \times  \frac{1}{10}(1, 3)) } $ \\

 \oddrow $\rownumber$ &$X_{36}\subset \PP(1, 7, 12, 18) $& 2 & $2/21$ & $1$ &
$11$ & $9$ & $ { \frac{1}{6}(1, 1), \frac{1}{7}(1, 3) } $ \\\hline\hline

 \evnrow $\rownumber$ &$X_{6}\subset \PP(1^2, 2, 5) $& 3 & $27/5$ & $6$ & $5$ &
$3$ & $ { \frac{1}{5}(1, 2) } $ \\

 \oddrow $\rownumber$ &$X_{8}\subset \PP(1^2,2,7 ) $& 3 & $36/7$ & $6$ & $6$
& $4$ & $ { \frac{1}{7}(1, 2) } $ \\

 \evnrow $\rownumber$ &$X_{10}\subset \PP(1, 2, 5^2) $& 3 & $9/5$ & $2$ & $7$
& $5$ & $ { 2 \times \frac{1}{5}(1, 2) } $ \\

 \oddrow $\rownumber$ &$X_{12}\subset \PP(1, 2, 5, 7) $& 3 & $54/35$ & $2$ & $8$
& $6$ & $ { \frac{1}{5}(1, 2), \frac{1}{7}(1, 2) } $ \\

 \evnrow $\rownumber$ &$X_{14}\subset \PP(1, 2, 7^2) $& 3 & $9/7$ & $2$ & $9$ &
$7$ & $ { 2 \times  \frac{1}{7}(1, 2) }$ \\

 \oddrow $\rownumber$ &$X_{15}\subset \PP(1, 5^2, 7) $& 3 & $27/35$ & $1$ & $9$
& $3$ & $ { 3 \times  \frac{1}{5}(1, 2), \frac{1}{7}(1, 1) } $\\\hline\hline

 \evnrow $\rownumber$ &$X_{6}\subset \PP(1^2, 3, 5) $& 4 & $32/5$ & $7$ & $4$ &
$2$ & $ { \frac{1}{5}(1, 2) } $ \\

 \oddrow $\rownumber$ &$X_{10}\subset \PP(1, 3, 5^2) $& 4 & $32/15$ & $2$ & $6$
& $4$ & $ { \frac{1}{3}(1, 1), 2 \times  \frac{1}{5}(1, 2) } $ \\

 \evnrow $\rownumber$ &$X_{12}\subset \PP(1, 3, 5, 7) $& 4 & $64/35$ & $2$ & $6$
& $4$ & $ { \frac{1}{5}(1, 2), \frac{1}{7}(1, 3) } $ \\

 \oddrow $\rownumber$ &$X_{15}\subset \PP(1, 3, 5, 10) $& 4 & $8/5$ & $2$ & $7$
& $5$ & $ { \frac{1}{5}(1, 2), \frac{1}{10}(1, 3) } $ \\

 \evnrow $\rownumber$ &$X_{15}\subset \PP(3, 5^2, 6) $& 4 & $8/15$ & $0$ & $7$ &
$5$ & $ { 2 \times  \frac{1}{3}(1, 1), 3 \times  \frac{1}{5}(1, 2), \frac{1}{6}(1, 1) } $ \\

 \oddrow $\rownumber$ &$X_{21}\subset \PP(1, 7^2, 10) $& 4 & $24/35$ & $1$ & $9$
& $7$ & $ { 3 \times  \frac{1}{7}(1, 3), \frac{1}{10}(1, 1) } $ \\\hline\hline

 \evnrow $\rownumber$ &$X_{8}\subset \PP(1, 2, 3, 7) $& 5 & $100/21$ & $5$ & $4$
& $2$ & $ { \frac{1}{3}(1, 1), \frac{1}{7}(1, 3) } $ \\

 \oddrow $\rownumber$ &$X_{8}\subset \PP(1^2, 4, 7) $& 5 & $50/7$ & $8$ & $4$ &
$2$ & $ { \frac{1}{7}(1, 2) } $ \\

 \evnrow $\rownumber$ &$X_{12}\subset \PP(1, 3, 4, 9) $& 5 & $25/9$ & $3$ & $5$
& $3$ & $ { \frac{1}{3}(1, 1), \frac{1}{9}(1, 4) } $ \\

 \oddrow $\rownumber$ &$X_{14}\subset \PP(2, 3, 7^2) $& 5 & $25/21$ & $1$ & $5$
& $3$ & $ { \frac{1}{3}(1, 1), 2 \times  \frac{1}{7}(1, 3) } $ \\

 \evnrow $\rownumber$ &$X_{16}\subset \PP(1, 4, 7, 9) $& 5 & $100/63$ & $2$ &
$6$ & $4$ & $ { \frac{1}{7}(1, 2), \frac{1}{9}(1, 4) } $ \\

 \oddrow $\rownumber$ &$X_{21}\subset \PP(3, 7^2, 9) $& 5 & $25/63$ & $0$ & $7$
& $5$ & $ { 2 \times  \frac{1}{3}(1, 1), 3 \times  \frac{1}{7}(1, 3), \frac{1}{9}(1, 1) } $ \\\hline\hline

 \evnrow $\rownumber$ &$X_{15}\subset \PP(1, 5, 7, 8) $& 6 & $27/14$ & $2$ & $5$
& $3$ & $ { \frac{1}{7}(1, 3), \frac{1}{8}(1, 5) } $ \\

 \oddrow $\rownumber$ &$X_{16}\subset \PP(1, 5, 8^2) $& 6 & $9/5$ & $2$ & $6$ &
$4$ & $ { \frac{1}{5}(1, 1), 2 \times  \frac{1}{8}(1, 5) } $ \\\hline\hline

 \evnrow $\rownumber$ &$X_{10}\subset \PP(1, 2, 5, 9) $& 7 & $49/9$ & $6$ & $3$
& $1$ & $ { \frac{1}{9}(1, 4) } $ \\

 \oddrow $\rownumber$ &$X_{12}\subset \PP(2, 3, 5, 9) $& 7 & $98/45$ & $2$ & $4$
& $2$ & $ { \frac{1}{3}(1, 1), \frac{1}{5}(1, 2), \frac{1}{9}(1, 4) } $ \\\hline\hline

 \evnrow $\rownumber$ &$X_{8}\subset \PP(1, 3, 5, 7) $& 8 & $512/105$ & $5$ &
$4$ & $2$ & $ { \frac{1}{3}(1, 1), \frac{1}{5}(1, 2), \frac{1}{7}(1, 2) } $ \\

 \oddrow $\rownumber$ &$X_{14}\subset \PP(1, 5, 7, 9) $& 8 & $128/45$ & $3$ &
$4$ & $2$ & $ { \frac{1}{5}(1, 2), \frac{1}{9}(1, 4) } $ \\

 \evnrow $\rownumber$ &$X_{15}\subset \PP(1, 6, 7, 9) $& 8 & $160/63$ & $3$ &
$5$ & $3$ & $ { \frac{1}{6}(1, 1), \frac{1}{7}(1, 3), \frac{1}{9}(1, 4) } $ \\

 \oddrow $\rownumber$ &$X_{15}\subset \PP(1, 5, 7, 10) $& 8 & $96/35$ & $3$ &
$5$ & $3$ & $ { \frac{1}{5}(1, 2), \frac{1}{7}(1, 2), \frac{1}{10}(1, 3)) } $ \\

\end{longtable}
 \renewcommand*{\arraystretch}{1.5}
\begin{longtable}{>{\hspace{0.5em}}llcccccc<{\hspace{0.5em}}}
\caption{Codimension 2 Complete Intersections  } \label{tab-CI}\\
\toprule
\multicolumn{1}{c}{S.No}&\multicolumn{1}{c}{$X$ }&\multicolumn{1}{c}{$I$}&\multicolumn{1}{c}{$-K^2$}&\multicolumn{1}{c}{$h^0(-K)$}&\multicolumn{1}{c}{$e(X)$}&\multicolumn{1}{c}{$e_\orb(X)$}&\multicolumn{1}{c}{Basket $\mathcal{B}$}\\
\cmidrule(lr){1-1}\cmidrule(lr){2-2}\cmidrule(lr){3-5}\cmidrule(lr){6-7}\cmidrule(lr){8-8}
\endfirsthead
\multicolumn{7}{l}{\vspace{-0.25em}\scriptsize\emph{\tablename\ \thetable{} continued from previous page}}\\
\midrule
\endhead
\multicolumn{7}{r}{\scriptsize\emph{Continued on next page}}\\
\endfoot
\bottomrule
\endlastfoot
\evnrow $\rownumber$ &$X_{2, 2}\subset \PP(1^5) $& 1 & $4$ & $5$& $8$ & $8$ & $ {
} $ \\ 

  \oddrow $\rownumber$ &$X_{4^2}\subset \PP(1^2, 2^2, 3) $& 1 & $4/3$ & $2$ & 
$10$ & $28/3$ & $ { \frac{1}{3}(1, 1) } $ \cite{dp-CH}\\ 

 \evnrow $\rownumber$ &$X_{4, 6}\subset \PP(1, 2^2, 3^2) $& 1 & $2/3$ & $1$ & 
$10$ & $26/3$ & $ { 2 \times  \frac{1}{3}(1, 1) } $\cite{dp-CH} \\ 

 \oddrow $\rownumber$ &$X_{6^2}\subset \PP(1^2, 3^2, 5) $& 1 & $4/5$ & $2$ & 
$12$ & $56/5$ & $ { \frac{1}{5}(1, 1) } $  \cite{Fano-CK}\\ 

 \evnrow $\rownumber$ &$X_{6^2}\subset \PP(2^2, 3^3) $& 1 & $1/3$ & $0$ & $9$ & 
$19/3$ & $ { 4 \times  \frac{1}{3}(1, 1) } $ \cite{dp-CH}\\ 

 \oddrow $\rownumber$ &$X_{6, 7}\subset \PP(1, 2, 3^2, 5) $& 1 & $7/15$ & $1$ & 
$11$ & $133/15$ & $ { 2 \times  \frac{1}{3}(1, 1), \frac{1}{5}(1, 1) } $ \\ 

 \evnrow $\rownumber$ &$X_{6, 8}\subset \PP(1, 2, 3, 4, 5) $& 1 & $2/5$ & $1$ & 
$10$ & $46/5$ & $ { \frac{1}{5}(1, 2) } $ \\ 

 \oddrow $\rownumber$ &$X_{8^2}\subset \PP(1^2, 4^2, 7) $& 1 & $4/7$ & $2$ & 
$14$ & $92/7$ & $ { \frac{1}{7}(1, 1) } $ \\ 

 \evnrow $\rownumber$ &$X_{6, 10}\subset \PP(1, 3^2, 5^2) $& 1 & $4/15$ & $1$ & 
$12$ & $136/15$ & $ { 2 \times  \frac{1}{3}(1, 1), 2 \times  \frac{1}{5}(1, 1) } $ \\ 

 \oddrow $\rownumber$ &$X_{8, 10}\subset \PP(2, 3, 4, 5^2) $& 1 & $2/15$ & $0$ &
$8$ & $86/15$ & $ { \frac{1}{3}(1, 1), 2 \times  \frac{1}{5}(1, 2) } $ \\ 

 \evnrow $\rownumber$ &$X_{9, 10}\subset \PP(2, 3^2, 5, 7) $& 1 & $1/7$ & $0$ & 
$9$ & $43/7$ & $ { 3 \times  \frac{1}{3}(1, 1), \frac{1}{7}(1, 2) } $ \\ 

 \oddrow $\rownumber$ &$X_{10^2}\subset \PP(1^2, 5^2, 9) $& 1 & $4/9$ & $2$ & 
$16$ & $136/9$ & $ { \frac{1}{9}(1, 1) } $ \\ 

 \evnrow $\rownumber$ &$X_{10, 11}\subset \PP(1, 2, 5^2, 9) $& 1 & $11/45$ & $1$
& $13$ & $473/45$ & $ { 2 \times  \frac{1}{5}(1, 2), \frac{1}{9}(1, 1) } $ \\ 

 \oddrow $\rownumber$ &$X_{10, 12}\subset \PP(3^2, 5^2, 7) $& 1 & $8/105$ & $0$ 
& $10$ & $512/105$ & $ { 4 \times  \frac{1}{3}(1, 1), 2 \times  \frac{1}{5}(1, 1), \frac{1}{7}(1, 2) } $ \\ 

 \evnrow $\rownumber$ &$X_{10, 12}\subset \PP(2, 3, 5, 6, 7) $& 1 & $2/21$ & $0$
& $8$ & $122/21$ & $ { 2 \times  \frac{1}{3}(1, 1), \frac{1}{7}(1, 3) } $ \\ \hline\hline

 \oddrow $\rownumber$ &$X_{6, 8}\subset \PP(1, 3^2, 4, 5) $& 2 & $16/15$ & $1$ &
$8$ & $88/15$ & $ { 2 \times  \frac{1}{3}(1, 1), \frac{1}{5}(1, 2) } $ \\ 

 \evnrow $\rownumber$ &$X_{8, 10}\subset \PP(1, 3, 4, 5, 7) $& 2 & $16/21$ & $1$
& $8$ & $136/21$ & $ { \frac{1}{3}(1, 1), \frac{1}{7}(1, 3) } $ \\ 

 \oddrow $\rownumber$ &$X_{8, 12}\subset \PP(1, 3, 5, 6, 7) $& 2 & $64/105$ & 
$1$ & $10$ & $736/105$ & $ { 2 \times  \frac{1}{3}(1, 1), \frac{1}{5}(1, 1), \frac{1}{7}(1, 2) } $ \\ 

 \evnrow $\rownumber$ &$X_{10, 12}\subset \PP(3, 4, 5^2, 7) $& 2 & $8/35$ & $0$ 
& $6$ & $124/35$ & $ { 2 \times  \frac{1}{5}(1, 2), \frac{1}{7}(1, 3)} $ \\ 

 \oddrow $\rownumber$ &$X_{9, 14}\subset \PP(1, 3, 6, 7, 8) $& 2 & $1/2$ & $1$ &
$10$ & $61/8$ & $ { \frac{1}{3}(1, 1), \frac{1}{6}(1, 1), \frac{1}{8}(1, 5) } $ \\ 

 \evnrow $\rownumber$ &$X_{12, 14}\subset \PP(3, 4, 5, 7, 9) $& 2 & $8/45$ & $0$
& $6$ & $164/45$ & $ { \frac{1}{3}(1, 1), \frac{1}{5}(1, 2), \frac{1}{9}(1, 4) } $ \\ 

 \oddrow $\rownumber$ &$X_{14, 15}\subset \PP(3, 6, 7^2, 8) $& 2 & $5/42$ & $0$ 
& $8$ & $545/168$ & $ { 2 \times  \frac{1}{3}(1, 1), \frac{1}{6}(1, 1), 2 \times  \frac{1}{7}(1, 2), \frac{1}{8}(1, 5) } $
\\ 

 \evnrow $\rownumber$ &$X_{11, 18}\subset \PP(1, 3, 8, 9, 10) $& 2 & $11/30$ & 
$1$ & $12$ & $1067/120$ & $ { 2 \times  \frac{1}{3}(1, 1), \frac{1}{8}(1, 1), \frac{1}{10}(1, 3) } $ \\ \hline\hline

 \oddrow $\rownumber$ &$X_{12, 14}\subset \PP(4,5,6,7^{2}) $& 3 & $9/35$ & 
$0$ & $5$ & $87/35$ & $ {  \frac{1}{5}(1, 2),2\times  \frac{1}{7}(1, 3) } $ \\ \hline\hline

 \evnrow $\rownumber$ &$X_{10, 12}\subset \PP(3, 5^2, 6, 7) $& 4 & $64/105$ & 
$0$ & $6$ & $232/105$ & $ { 2 \times  \frac{1}{3}(1, 1), 2 \times  \frac{1}{5}(1, 2), \frac{1}{7}(1, 2) } $ \\

 \end{longtable}
\renewcommand*{\arraystretch}{.9}
\begin{longtable}{>{\hspace{0.5em}}llccccr<{\hspace{0.5em}}}
\caption{Codimension 3 Pfaffians\  } \label{tab-Pf}\\
\toprule
\multicolumn{1}{c}{S.No}&\multicolumn{1}{c}{$X$ }&\multicolumn{1}{c}{$I$}&\multicolumn{1}{c}{$-K^2$}&\multicolumn{1}{c}{$h^0(-K)$}&\multicolumn{1}{c}{Basket $\mathcal{B}$}&\multicolumn{1}{c}{Weight Matrix }\\
\cmidrule(lr){1-1}\cmidrule(lr){2-2}\cmidrule(lr){3-5}\cmidrule(lr){6-6}\cmidrule(lr){7-7}
\endfirsthead
\multicolumn{7}{l}{\vspace{-0.25em}\scriptsize\emph{\tablename\ \thetable{} continued from previous page}}\\
\midrule
\endhead
\multicolumn{7}{r}{\scriptsize\emph{Continued on next page}}\\
\endfoot
\bottomrule
\endlastfoot

\oddrow $\rownumber$ &$\begin{array}{@{}l@{}} X_{2, 2, 2, 2, 2}\\ \quad\subset 
\PP(1^6)\end{array}$& 1 & $5$ & $6$&  & $\begin{matrix} 1&1&1&1\\ &1&1&1\\ &&1&1\\ &&&1 \end{matrix}$ \\ 

  \evnrow $\rownumber$ &$\begin{array}{@{}l@{}} X_{3, 3, 4, 4, 4}\\ \quad\subset 
\PP(1^3, 2^2, 3)\end{array} $& 1 & $7/3$ & $3$& $ { \frac13(1, 1) } $~\cite{dp-CH} &$\begin{matrix} 1&1&2&2\\ &1&2&2\\ &&2&2\\ &&&3 
\end{matrix}$ \\ 

  \oddrow $\rownumber$ &$\begin{array}{@{}l@{}} X_{4, 4, 6, 6, 6}\\ \quad\subset 
\PP(1^3, 3^2, 5)\end{array} $& 1 & $9/5$ & $3$& $ { \frac15(1, 1) } $ \cite{Fano-CK} & 
$\begin{matrix} 1&1&3&3\\ &1&3&3\\ &&3&3\\ &&&5 
\end{matrix}$ \\ 

  \evnrow $\rownumber$ &$\begin{array}{@{}l@{}} X_{4, 5, 6, 6, 7}\\ \quad\subset 
\PP(1^2, 2, 3^2, 5)\end{array} $& 1 & $17/15$ & $2$& $ { \frac13(1, 1), \frac15(1, 1) } 
$ &$\begin{matrix} 1&1&2&3\\ &2&3&4\\ &&3&4\\ &&&5 
\end{matrix}$ \\ 

  \oddrow $\rownumber$ &$\begin{array}{@{}l@{}} X_{5, 5, 8, 8, 8}\\ \quad\subset 
\PP(1^3, 4^2, 7)\end{array} $& 1 & $11/7$ & $3$& $ { \frac{1}{7}(1, 1) } $ & 
$\begin{matrix} 1&1&4&4\\ &1&4&4\\ &&4&4\\ &&&7 
\end{matrix}$ \\ 


  \evnrow $\rownumber$ &$\begin{array}{@{}l@{}} X_{6, 7, 8, 9, 10}\\ \quad\subset 
\PP(1, 2, 3^2, 5, 7)\end{array} $& 1 & $10/21$ & $1$& $ { \frac13(1, 1), \frac17(1, 4) }
$ &$\begin{matrix} 1&2&3&4\\ &3&4&5\\ &&5&6\\ &&&7 
\end{matrix}$ \\ 

  \oddrow $\rownumber$ &$\begin{array}{@{}l@{}} X_{6, 6, 10, 10, 10}\\ \quad\subset 
\PP(1^3, 5^2, 9)\end{array} $& 1 & $13/9$ & $3$& $ {\frac19(1, 1) } $ & 
$\begin{matrix} 1&1&5&5\\ &1&5&5\\ &&5&5\\ &&&9 
\end{matrix} $\\ 

  \evnrow $\rownumber$ &$\begin{array}{@{}l@{}} X_{7, 8, 8, 9, 10}\\ \quad\subset 
\PP(2, 3^2, 4, 5^2)\end{array} $& 1 & $1/5$ & $0$& $ { 3 \times \frac13(1, 1), \frac15(1, 2),
\frac15(1, 1) } $ & $\begin{matrix} 2&3&3&4\\ &4&4&5\\ 
&&5&6\\ &&&6 \end{matrix}$ \\ 

  \oddrow $\rownumber$ &$\begin{array}{@{}l@{}} X_{6, 7, 10, 10, 11}\\ \quad\subset 
\PP(1^2, 2, 5^2, 9)\end{array} $& 1 & $38/45$ & $2$& $ { \frac15(1, 2),\frac19(1, 1) } 
$ & $\begin{matrix} 1&1&4&5\\ &2&5&6\\ &&5&6\\ &&&9 
\end{matrix}$ \\ 

  \evnrow $\rownumber$ &$\begin{array}{@{}l@{}} X_{6, 8, 10, 10, 12}\\ \quad\subset 
\PP(1, 3^2, 5^2, 7)\end{array} $& 1 & $29/105$ & $1$& $ { \frac13(1, 1), \frac15(1, 1), 
\frac17(1, 4) } $ & $\begin{matrix} 1&1&3&5\\ 
&3&5&7\\ &&5&7\\ &&&9 \end{matrix} $\\ 

  \oddrow $\rownumber$ &$\begin{array}{@{}l@{}} X_{10, 10, 12, 12, 14}\\ \quad\subset 
\PP(3^2, 5^2, 7^2)\end{array} $& 1 & $3/35$ & $0$& $ { 3 \times \frac13(1, 1), \frac15(1, 1),
2 \times \frac17(1, 4) } $ & $\begin{matrix} 3&3&5&5\\ 
&5&7&7\\ &&7&7\\ &&&9 \end{matrix}$ \\ 

  \evnrow $\rownumber$ &$\begin{array}{@{}l@{}} X_{11, 12, 12, 15, 16}\\ \quad\subset 
\PP(2, 5^2, 6, 7, 9)\end{array} $& 1 & $23/315$ & $0$& $ { 3 \times \frac15(1, 2), \frac{1}{7}(1, 3),\frac19(1, 1) } $ & $\begin{matrix} 2&5&5&6\\ 
&6&6&7\\ &&9&10\\ &&&10 \end{matrix}$ \\ 

  \oddrow $\rownumber$ &$\begin{array}{@{}l@{}} X_{4, 7, 8, 8, 9}\\ \quad\subset 
\PP(1^2, 2, 3, 6, 7)\end{array} $& 2 & $22/7$ & $4$& $ {\frac13(1, 1), \frac16(1, 1), 
\frac{1}{7}(1, 2) } $ &$\begin{matrix} 1&1&2&5\\ &2&3&6\\ &&3&6\\ 
&&&7 \end{matrix}$ \\ 

  \evnrow $\rownumber$ &$\begin{array}{@{}l@{}} X_{4, 8, 9, 9, 10}\\ \quad\subset 
\PP(1^2, 2, 3, 7, 8)\end{array} $& 2 & $43/14$ & $4$& $ { \frac17(1, 1), \frac{1}{8}(1, 5) }
$ &$\begin{matrix} 1&1&2&6\\ &2&3&7\\ &&3&7\\ &&&8 
\end{matrix}$ \\ 

\oddrow $\rownumber$ &$\begin{array}{@{}l@{}} X_{4, 10, 11, 11, 12}\\ \quad\subset 
\PP(1^2, 2, 3, 9, 10)\end{array} $& 2 & $134/45$ & $4$& $ {\frac13(1, 1), \frac{1}{9}(1, 1), \frac{1}{10}(1, 3) } $ &$\begin{matrix} 1&1&2&8\\ &2&3&9\\ 
&&3&9\\ &&&10 \end{matrix}$ \\ 

\evnrow $\rownumber$ &$\begin{array}{@{}l@{}} X_{8, 9, 12, 13, 14}\\ \quad\subset 
\PP(1, 3, 5, 6, 7, 8)\end{array} $& 2 & $19/30$ & $1$& $ {\frac13(1, 1), \frac15(1, 1),
\frac{1}{8}(1, 5) } $ &$\begin{matrix} 1&2&5&6\\ &3&6&7\\ &&7&8\\ 
&&&11 \end{matrix}$ \\ 

\oddrow $\rownumber$ &$\begin{array}{@{}l@{}} X_{12, 12, 14, 15, 15}\\ \quad\subset 
\PP(4, 5^2, 7^2, 8)\end{array} $& 2 & $11/70$ & $0$& $ { 2 \times \frac{1}{5}(1, 2), 2 \times 
\frac{1}{7}(1, 3), \frac{1}{8}(1, 1) } $ &$\begin{matrix} 4&5&7&7\\ &5&7&7\\ 
&&8&8\\ &&&10 \end{matrix}$ \\ 


  \evnrow $\rownumber$ &$\begin{array}{@{}l@{}} X_{14, 14, 15, 15, 16}\\ \quad\subset 
\PP(3, 6, 7^2, 8^2)\end{array} $& 2 & $1/7$ & $0$& $ {\frac13(1, 1), \frac16(1, 1), 
\frac{1}{7}(1, 2), 2 \times \frac{1}{8}(1, 5) } $ & $\begin{matrix} 
6&6&7&7\\ &7&8&8\\ &&8&8\\ &&&9 \end{matrix}$ \\ 

  \oddrow $\rownumber$ &$\begin{array}{@{}l@{}} X_{16, 17, 17, 18, 18}\\ \quad\subset 
\PP(3, 7, 8^2, 9, 10)\end{array} $& 2 & $11/105$ & $0$& $ {\frac13(1, 1), \frac{1}{7}(1, 1), 2 \times \frac{1}{8}(1, 5), \frac{1}{10}(1, 3) } $ & $\begin{matrix} 7&8&8&9\\ &8&8&9\\ &&9&10\\ &&&10 \end{matrix}$ \\ \hline\hline
\evnrow $\rownumber$ &$\begin{array}{@{}l@{}} X_{6^,6,8,8, 10}\\ \quad\subset 
\PP(1^2, 3, 5^2, 7)\end{array} $& 3 & $153/35$ & $4$& $ { 2 \times \frac{1}{5}(1, 2), 
\frac{1}{7}(1, 1) } $ &$\begin{matrix} 1&1&3&3\\ &3&5&5\\ &&5&5\\ 
&&&7 \end{matrix} $\\ \hline\hline

  \oddrow $\rownumber$ &$\begin{array}{@{}l@{}} X_{8, 8, 11, 11, 14}\\ \quad\subset 
\PP(1^2, 4, 7^2, 10)\end{array} $& 4 & $184/35$ & $6$& $ { 2 \times \frac{1}{7}(1, 3), 
\frac{1}{10}(1, 1) } $ &$\begin{matrix} 1&1&4&4\\ &4&7&7\\ &&7&7\\ 
&&&10 \end{matrix}$ \\ 
 
\end{longtable}

\renewcommand*{\arraystretch}{.9}
\begin{longtable}{>{\hspace{0.5em}}llccccc<{\hspace{0.5em}}}
\caption{Codimension 4 \(\PxP\)\  } \label{tab-P2xP2}\\
\toprule
\multicolumn{1}{c}{S. No}&\multicolumn{1}{c}{$X$ }&\multicolumn{1}{c}{$I$}&\multicolumn{1}{c}{$-K^2$}&\multicolumn{1}{c}{$h^0(-K)$}&\multicolumn{1}{c}{Basket $\mathcal{B}$}&\multicolumn{1}{c}{Weight Matrix }\\
\cmidrule(lr){1-1}\cmidrule(lr){2-2}\cmidrule(lr){3-5}\cmidrule(lr){6-6}\cmidrule(lr){7-7}
\endfirsthead
\multicolumn{7}{l}{\vspace{-0.25em}\scriptsize\emph{\tablename\ \thetable{} continued from previous page}}\\
\midrule
\endhead
\multicolumn{7}{r}{\scriptsize\emph{Continued on next page}}\\
\endfoot
\bottomrule
\endlastfoot

\evnrow $\rownumber$ &$\begin{array}{@{}l@{}} X_{2^{9}}\\ \quad\subset \PP(1^7)\end{array} $& 1 & $6$ & $7$& $ { } $ & 
$\begin{matrix} 1&1&1\\ 1&1&1\\ 1&1&1 \end{matrix}$ \\ 

  \oddrow $\rownumber$ &$\begin{array}{@{}l@{}} X_{2,3^{4},4^{4}}\\ \quad\subset \PP(1^4, 2^2, 3)\end{array} $& 1 & $10/3$ & $4$&
   $  \frac13(1, 1) $
  \cite{dp-CH} & $\begin{matrix} 1&1&2\\ 1&1&2\\ 2&2&3 \end{matrix}$ \\ 

  \evnrow $\rownumber$ &$\begin{array}{@{}l@{}} X_{2, 4^{4}, 6^{4}}\\ \quad\subset \PP(1^4, 3^2, 5)\end{array} $& 1 & $14/5$ & $4$& $ { \frac15(1, 1)
} $  \cite{Fano-CK}& $\begin{matrix} 1&1&3\\ 1&1&3\\ 3&3&5 \end{matrix}$ \\ 

  \oddrow $\rownumber$ &$\begin{array}{@{}l@{}} X_{2, 5^{4}, 8^{4}}\\ \quad\subset \PP(1^4, 4^2, 7)\end{array} $& 1 & $18/7$ & $4$& $ { \frac{1}{7}(1, 1)
} $ & $\begin{matrix} 1&1&4\\ 1&1&4\\ 4&4&7 \end{matrix}$ \\ 


  \evnrow $\rownumber$ &$\begin{array}{@{}l@{}} X_{4, 5^{2}, 6^{3}, 7^{2}, 
8}\\ \quad\subset \PP(1, 2^2, 3^3, 5)\end{array} $& 1 & $4/5$ & $1$& $ { 3 \times 
\frac13(1, 1), \frac15(1, 1) } $ & $\begin{matrix} 1&2&3\\ 2&3&4\\ 3&4&5 \end{matrix}$ 
\\ 

  \oddrow $\rownumber$ &$\begin{array}{@{}l@{}} X_{2, 6^{4}, 10^{4}}\\ \quad\subset \PP(1^4, 5^2, 9)\end{array} $& 1 & $22/9$ & $4$& $ { \frac{1}{9}(1, 1) } $ & $\begin{matrix} 1&1&5\\ 1&1&5\\ 5&5&9 \end{matrix}$ \\ 

  \evnrow $\rownumber$ &$\begin{array}{@{}l@{}} X_{5, 6^{2}, 7^{2}, 8^{2}, 9, 
10}\\ \quad\subset \PP(1, 2, 3^2, 4, 5^2)\end{array} $& 1 & $8/15$ & $1$& $ { 
\frac13(1, 1), \frac15(1, 2), \frac15(1, 1) } $ & $\begin{matrix} 1&2&3\\ 3&4&5\\ 4&5&6 
\end{matrix}$ \\ 


  \oddrow $\rownumber$ &$\begin{array}{@{}l@{}} X_{4, 7^{2}, 8^{2}, 10, 11^{2}, 
12}\\ \quad\subset \PP(1, 2^2, 3, 5^2, 9)\end{array} $& 1 & $26/45$ & $1$& $ { 
\frac13(1, 1), 2 \times \frac15(1, 2),\frac19(1, 1) } $ & $\begin{matrix} 1&2&5\\ 2&3&6\\ 5&6&9 
\end{matrix}$ \\ 


  \evnrow $\rownumber$ &$\begin{array}{@{}l@{}} X_{6, 8^{2}, 10^{3}, 12^{2}, 
14}\\ \quad\subset \PP(1, 3^2, 5^2, 7^2)\end{array} $& 1 & $2/7$ & $1$& $ { 2 \times 
\frac17(1, 4) } $ & $\begin{matrix} 1&3&5\\ 3&5&7\\ 5&7&9 \end{matrix}$ \\ 

  \oddrow $\rownumber$ &$\begin{array}{@{}l@{}} X_{7, 8, 10, 11^{2}, 12^{2}, 15, 
16}\\ \quad\subset \PP(1, 2, 5^2, 6, 7, 9)\end{array} $& 1 & $86/315$ & $1$& $ {
\frac15(1, 2), \frac{1}{7}(1, 3),\frac19(1, 1) } $ & $\begin{matrix} 1&2&5\\ 5&6&9\\ 6&7&10 
\end{matrix}$ \\ 

  \evnrow $\rownumber$ &$\begin{array}{@{}l@{}} X_{10, 11^{2}, 12^{3}, 13^{2}, 
14}\\ \quad\subset \PP(3, 4, 5^2, 6, 7^2)\end{array} $& 1 & $3/35$ & $0$& $ { 2 
\times \frac15(1, 2), 2 \times \frac17(1, 4) } $ & $\begin{matrix} 4&5&6\\ 5&6&7\\ 6&7&8 
\end{matrix}$ \\ 

  \oddrow $\rownumber$ &$\begin{array}{@{}l@{}} X_{8, 9, 11, 12^{3}, 13, 15, 
16}\\ \quad\subset \PP(2, 3, 5^2, 6, 7, 9)\end{array} $& 1 & $38/315$ & $0$& $ {
\frac13(1, 1), 3 \times \frac15(1, 2), \frac17(1, 4),\frac19(1, 1) } $ & $\begin{matrix} 2&3&6\\ 
5&6&9\\ 6&7&10 \end{matrix}$ \\ 


  \evnrow $\rownumber$ &$\begin{array}{@{}l@{}} X_{4, 8^{2}, 9^{2}, 12, 13^{2}, 
14}\\ \quad\subset \PP(1, 2, 3, 6^2, 7^2)\end{array} $& 2 & $20/21$ & $2$& $ { 
\frac13(1, 1), 2 \times \frac16(1, 1), 2 \times \frac{1}{7}(1, 2) } $ & $\begin{matrix} 1&2&6\\ 2&3&7\\ 
6&7&11 \end{matrix}$ \\ 

\oddrow $\rownumber$ &$\begin{array}{@{}l@{}} X_{4, 8, 9^{2}, 10, 13, 14^{2}, 
15}\\ \quad\subset \PP(1, 2, 3, 6, 7^2, 8)\end{array} $& 2 & $37/42$ & $2$& $ { 
\frac16(1, 1), \frac17(1, 1), \frac{1}{7}(1, 2), \frac{1}{8}(1, 5) } $ & $\begin{matrix} 1&2&6\\ 2&3&7\\
7&8&12 \end{matrix}$ \\ 

\evnrow $\rownumber$ &$\begin{array}{@{}l@{}} X_{4, 8, 9, 11, 12, 15, 16^{2}, 
17}\\ \quad\subset \PP(1, 2, 3, 6, 7, 9, 10)\end{array} $& 2 & $248/315$ & $2$& 
$ {\frac13(1, 1), \frac16(1, 1), \frac{1}{7}(1, 2), \frac19(1, 1), \frac{1}{10}(1, 3) } $ & 
$\begin{matrix} 1&2&6\\ 2&3&7\\ 9&10&14 \end{matrix}$ \\ 

  \oddrow $\rownumber$ &$\begin{array}{@{}l@{}} X_{4, 9, 10, 11, 12, 16, 17^{2}, 
18}\\ \quad\subset \PP(1, 2, 3, 7, 8, 9, 10)\end{array} $& 2 & $451/630$ & $2$& 
$ { \frac17(1, 1), \frac{1}{8}(1, 5), \frac19(1, 1), \frac{1}{10}(1, 3) } $ & $\begin{matrix} 1&2&7\\ 
2&3&8\\ 9&10&15 \end{matrix}$ \\ 

  \evnrow $\rownumber$ &$\begin{array}{@{}l@{}} X_{4, 11^{2}, 12^{2}, 18, 19^{2}, 
20}\\ \quad\subset \PP(1, 2, 3, 9^2, 10^2)\end{array} $& 2 & $28/45$ & $2$& $ { 
\frac13(1, 1), 2 \times \frac19(1, 1), 2 \times \frac{1}{10}(1, 3) } $ & $\begin{matrix} 1&2&9\\ 2&3&10\\ 
9&10&17 \end{matrix}$ \\ 


  \oddrow $\rownumber$ &$\begin{array}{@{}l@{}} X_{14, 15^{2}, 16^{3}, 17^{2}, 
18}\\ \quad\subset \PP(3, 6, 7^2, 8, 9, 10)\end{array} $& 2 & $16/105$ & $0$& $ 
{ 3 \times\frac13(1, 1), \frac16(1, 1), 2 \times \frac{1}{7}(1, 2), \frac{1}{10}(1, 3) } $ & $\begin{matrix} 
6&7&8\\ 7&8&9\\ 8&9&10 \end{matrix}$ \\ \hline\hline

  \evnrow $\rownumber$ &$\begin{array}{@{}l@{}} X_{6, 8^{2}, 10^{3}, 12^{2}, 
14}\\ \quad\subset \PP(1, 3^2, 5, 7^2, 9)\end{array} $& 5 & $250/63$ & $4$& $ { 
2 \times\frac13(1, 1), 2 \times \frac{1}{7}(1, 3), \frac{1}{9}(1, 1) } $ & $\begin{matrix} 1&3&5\\ 3&5&7\\ 
5&7&9 \end{matrix}$ \\ \hline\hline

  \oddrow $\rownumber$ &$\begin{array}{@{}l@{}} X_{14, 15^{2}, 16^{3}, 17^{2}, 
18}\\ \quad\subset \PP(6, 7^2, 8, 9^2, 10)\end{array} $& 8 & $256/315$ & $1$& $ 
{ \frac16(1, 1), 2 \times \frac{1}{7}(1, 3), 2 \times \frac{1}{9}(1, 4), \frac{1}{10}(1, 1) } $ & $\begin{matrix} 
6&7&8\\ 7&8&9\\ 8&9&10 \end{matrix}$ \\ 

\end{longtable}


\end{document}